\documentclass{article} 
   \usepackage{amssymb}  
\newcommand{\Aa}[1]{\overline{#1}}
\newcommand{\Ae}[1]{\overline{\overline{#1}}}
\newcommand{\Asp}[1]{\overline{\overline{\overline{#1}}}}
\newcommand{\Gq}{\mbox{$ \Ae{G} $}}
   
\newcommand{\Tq}{\mbox{$ \Ae{T} $}}
\newcommand{\g}{\mbox{$\bf g$}}

\newcommand{\h}{\mbox{\textbf{h}}}
\newcommand{\n}{\mbox{\textbf{n}}}
\newcommand{\np}{\mbox{$\textbf{n}^+$}}
\newcommand{\nm}{\mbox{$\textbf{n}^-$}}
\newcommand{\al}{\alpha}
\newcommand{\eps}{\epsilon}
\newcommand{\la}{\lambda}
\newcommand{\La}{\Lambda}
\newcommand{\de}{\delta}
\newcommand{\Th}{\Theta}
\newcommand{\gt}{\theta}
\newcommand{\mb}{\mbox}

\newcommand{\Mklz}[2]{\left\{\left.\;#1\;\right|\; #2\;\right\}}
\newcommand{\W}{\mbox{$\Delta$}}
\newcommand{\nW}{\mbox{$\Delta^-$}}
\newcommand{\pW}{\mbox{$\Delta^+$}}
\newcommand{\rW}{\mbox{$\Delta_{re}$}}

\newcommand{\prW}{\mbox{$\Delta_{re}^+$}}
\newcommand{\iW}{\mbox{$\Delta_{im}$}}

\newcommand{\Ht}{\mbox{ht}}
\newcommand{\F}{\mathbb{F}}

\newcommand{\C}{\mathbb{C}}
\newcommand{\N}{\mathbb{N}}
\newcommand{\Nn}{\mathbb{N}_0}
\newcommand{\Q}{\mathbb{Q}}
\newcommand{\Qp}{\mathbb{Q}^+}

\newcommand{\R}{\mathbb{R}}

\newcommand{\Z}{\mathbb{Z}}
\newcommand{\We}{\mbox{$\cal W$}}

\newcommand{\Cq}{\mbox{$\overline{C}$}}

\newcommand{\RkX}{\mbox{${\cal R}(X)$}}
\newcommand{\iB}[2]{\left(#1\mid#2\right)}
\newcommand{\iBl}{\left( \;\ \mid \;\ \right)}
\newcommand{\kB}[2]{\left\langle \left\langle #1\mid #2
                       \right\rangle \right\rangle}
\newcommand{\kkB}[2]{\langle \langle #1\mid #2\rangle \rangle}                       
\newcommand{\kBl}{\left\langle \left\langle \;\ \mid \;\ \right\rangle 
                       \right\rangle}
\newcommand{\grAdj}{\mbox{$gr$-$Adj$}}
\newcommand{\grEnd}{\mbox{$gr$-$End$}}
              
\newcommand{\Endo}{\mbox{$End$}\left(\,\bigoplus_{\La\in P^+}\,L(\La)\,
                    \right)}
\newcommand{\TD}{\mbox{$\widehat{T}$}}
\newcommand{\ND}{\mbox{$\widehat{N}$}}
\newcommand{\GD}{\mbox{$\widehat{G}$}}
\newcommand{\WeD}{\mbox{$\widehat{\We}$}}
\newcommand{\ve}[1]{\mbox{$\varepsilon\left(#1\right)$}} 
\newcommand{\FK}[1]{\mbox{$\F\,[#1]$}}
\newcommand{\CK}[1]{\mbox{${\mathbb C}\,[#1]$}}

\newcommand{\Proof}{\mbox{\bf Proof: }}
\newcommand{\End}{\mb{}\hfill\mb{$\square$}\\}
\newcommand{\Ende}{\mb{}\hfill\mb{$\square$}\vspace*{1ex}\\}
\newcommand{\Spm}{\mbox{Specm\,}}

\newcommand{\ti}{\tilde}
\newcommand{\res}[1]{\!\mid_{#1}}

\newcommand{\tr}{{\,\bf \diamond\,}}
\begin{document}
\newtheorem{Theorem}{Theorem}[section]
\newtheorem{Def}[Theorem]{Definition}
\newtheorem{Prop}[Theorem]{Proposition}
\newtheorem{Prop+Def}[Theorem]{Proposition+Definition}
\newtheorem{Cor}[Theorem]{Corollary}
\newtheorem{Rem}[Theorem]{Remark}
\newtheorem{Rems}[Theorem]{Remarks}
\newtheorem{Lemma}[Theorem]{Lemma}
%
%
\title{The $\F$-valued points of the algebra of strongly regular functions of a Kac-Moody group}
\author{Claus Mokler\\\\ Universit\"at Wuppertal, Fachbereich Mathematik\\  Gau\ss stra\ss e 20\\ D-42097 Wuppertal, Germany\vspace*{1ex}\\ 
          mokler@math.uni-wuppertal.de}
\date{}
\maketitle
\begin{abstract}\noindent
Let $G_m$ resp. $G_f$ be the minimal resp. formal Kac-Moody group, associated to a symmetrizable generalized Cartan matrix,
over a field $\F$ of characteristic 0. Let $\FK{G_m}$ be the algebra of strongly regular functions on $G_m$.\\
We denote by $\widehat{G_m}$ resp. $\widehat{G_f}$ certain monoid completions of $G_m$ resp. $G_f$, build by using the faces of 
the Tits cone.\\
We show that there is an action of $\widehat{G_f}\times\widehat{G_f}$ on the spectrum of $\F$-valued points of $\FK{G_m}$. As a 
$\widehat{G_f}\times\widehat{G_f}$-set it can be identified with a certain quotient of the $\widehat{G_f}\times\widehat{G_f}$-set 
$\widehat{G_f}\times\widehat{G_f}$, build by using $\widehat{G_m}$.\\
We prove a Birkhoff decomposition for the $\F$-valued points of $\FK{G_m}$.\\
We describe the stratification of the spectrum of $\F$-valued points of $\FK{G_m}$ in $G_f\times G_f$-orbits. We show that every 
orbit can be covered by suitably defined big cells.
\end{abstract}
\vspace*{1ex}
\section*{Introduction}
The minimal Kac-Moody group $G_m$, which V. Kac and D. Peterson associated in \cite{KP1} to a Kac-Moody algebra $\g$ over a 
field $\F$ of characteristic 0, is a group analogue of a semisimple simply connected algebraic group.\\
For a symmetrizable minimal Kac-Moody group, Kac and Peterson defined and investigated in  \cite{KP2} the algebra of strongly regular 
functions $\FK{G_m}$ on $G_m$. This algebra has many properties in common with the coordinate ring of a semisimple simply connected
algebraic group. It is an integrally closed domain, even a unique factorization domain. It admits a Peter and Weyl theorem, i.e.,
\begin{eqnarray*}
  \FK{G_m} &\cong & \bigoplus_{\La\in P^+} L^*(\La)\otimes L(\La)
\end{eqnarray*}
as $G_m\times G_m$-modules. But the following things, which hold in the non-classical case, are different:\\
1) Assigning to every element of $G_m$ its point evaluation, $G_m$ embeds in the set of $\F$-valued points of $\FK{G_m}$, which we 
denote by $\Spm\FK{G_m}$. But this map is not surjective.\\
2) There exists no comultiplication of $\FK{G_m}$, dual to the multiplication of $G_m$.  
The situation is not too bad, left and right multiplications with elements of $G_m$ induce comorphisms.
A more serious difference, the inverse map of $G_m$ does not induce a comorphism.\\
In particular there is no natural group structure, even no natural monoid structure on $\Spm\FK{G_m}$.\vspace*{0.5ex}\\
Kac and Peterson posed the problem to determine $\Spm\FK{G_m}$, or at least a certain part of it, \cite{KP2}.\\\\
The tensor category ${\cal O}_{adm}$ of admissible modules of $\cal O$ generalizes the category of finite dimensional 
representations of a semisimple Lie algebra, keeping the complete reducibility theorem.\\
In a way similar to the reconstruction process of the Tannaka-Krein duality, a suitable category of representations of a Lie algebra, together with 
a suitable category of duals, determines a monoid with coordinate ring. In some sense, this monoid is the biggest monoid acting 
reasonably on the representations. The coordinate ring is a coordinate ring of matrix coefficients.\\
We determined and investigated in \cite{M} the monoid $\widehat{G_m}$ corresponding to ${\cal O}_{adm}$ and its category of restricted duals. Equipped 
with its coordinate ring $\FK{\widehat{G_m}}$ of matrix coefficients, the monoid $\widehat{G_m}$ contains $G_m$ as Zariski open, dense unit group. It 
has similar properties as a reductive algebraic monoid. But it is a purely non-classical phenomenon, its classical analogue is a semisimple, simply 
connected algebraic group. For generalizing results of classical invariant theory, this monoid is more fundamental than the Kac-Moody group itself. For 
its history in connection with Slodowy and Peterson we refer to the introduction of \cite{M}.\\
The coordinate ring $\FK{\widehat{G_m}}$ is isomorphic to the algebra of strongly regular functions $\FK{G_m}$ by the
restriction map. Therefore the monoid $\widehat{G_m}$ embeds in $\Spm\FK{G_m}$, but also this map is not surjective.\vspace*{1ex}\\
To investigate the $\F$-valued points of $\FK{G_m}$, we define and investigate a monoid $\widehat{G_f}$, which is build in a similar way 
as $\widehat{G_m}$, but the minimal Kac-Moody group $G_m$ replaced by the formal Kac-Moody group $G_f$. In a subsequent paper we will prove 
that this monoid corresponds to ${\cal O}_{adm}$ and its category of full duals.\\
In this paper we obtain the following description of $\Spm\FK{G_m}$: 
We get, in a natural way, an action $\pi$ of $\widehat{G_f}\times \widehat{G_f}$ on $\FK{G_m}$ by homomorphisms of algebras. Therefore we 
also obtain a $\widehat{G_f}\times \widehat{G_f}$-action on the spectrum of $\F$-valued points of $\FK{G_m}$ from the right. Composing the 
evaluation map at the unit of $G_m$ with $\pi$, we get a $\widehat{G_f}\times \widehat{G_f}$-equivariant map 
\begin{eqnarray*}
  \tr:\,\widehat{G_f}\times \widehat{G_f} \,\;\to \;\, \Spm\FK{G_m}   & \quad,\quad &  (x,y) \mapsto x\tr y  \;\;.
\end{eqnarray*}
We show that this map factors to a $\widehat{G_f}\times \widehat{G_f}$-equivariant bijection with a quotient set of $\widehat{G_f}\times \widehat{G_f}$, 
which is obtained as follows: The Chevalley involution of $G_m$ extends to an involution $*$ of $\widehat{G_m}$. We factor 
$\widehat{G_f}\times \widehat{G_f}$ by the $\widehat{G_f}\times \widehat{G_f}$-equivariant equivalence relation generated by
\begin{eqnarray*}
  (x,zy) \;\sim\;(z^*x,y ) &\quad,\quad & x,y\in \widehat{G_f}\;,\; z\in \widehat{G_m}\;.
\end{eqnarray*}
Due to this description, we can use structural properties of $\widehat{G_m}$ and $\widehat{G_f}$ to prove properties of the 
$\widehat{G_f}\times \widehat{G_f}$-space $\Spm\FK{G}$.
In particular we prove the Birkhoff decomposition
\begin{eqnarray*}
  \Spm\FK{G_m} &=& \dot{\bigcup_{\hat{w}\in \widehat{\cal W}}} B_f\tr \hat{w} B_f\;\;.
\end{eqnarray*}
Here $B_f$ is the formal Borel group, and the Weyl monoid $\WeD$ is a certain monoid containing the Weyl group.\\ 
We determine the stratification of $\Spm\FK{G}$ in $G_f\times G_f$-orbits:
\begin{eqnarray*}
  \Spm\FK{G_m} &=& \dot{\bigcup_{\Th\,special}} G_f\tr e(R(\Th)) G_f\;\;.
\end{eqnarray*}
Here $\Mklz{e(R(\Th))}{\Th\mb{ special }}$ is a finite set of certain idempotents of $\widehat{G_m}$. We show that each orbit is locally closed and 
irreducible. We determine the closure relation of the orbits. We describe the big cell $B_f\tr e(R(\Th)) B_f$ of the orbit $G_f \tr e(R(\Th)) G_f$, and 
also the covering of $G_f \tr e(R(\Th)) G_f$ by this big cell. We give stratified transversal slices to the orbits.\vspace*{1ex}\\
The following work is in relation to these results: \vspace*{0.5ex}\\
Let $(M,\CK{M})$ be a connected reductive algebraic monoid. Denote by $G$ its reductive unit group. Let $T$ be a maximal torus of $G$. Let $B$ 
and $B^-$ be opposite Borel subgroups containing $T$.\\ 
Assigning to every element of $M$ its point evaluation, $M$ identifies with the $\C$-valued points of $\CK{M}$.\\
L. Renner showed in \cite{Re2}, that $M$ admits Bruhat decompositions. Due to the existence of a longest element of the Weyl group, these are 
equivalent to the Birkhoff decompositions. In particular $M=\dot{\bigcup}_{r\in {\cal R}} B^- r B$, where ${\cal R}$ is the Renner monoid.
Due to the work of M. Putcha and L. Renner in \cite{Pu2}, \cite{Re1} is the decomposition $M=\dot{\bigcup}_{e\in \Lambda} GeG$. The set $\Lambda$ 
is a certain set of idempotents, a cross-section lattice, which has been introduced by Putcha in \cite{Pu1}. 
Renner found and described in \cite{Re2} the big cell $B^-eB$ in $GeG$.\vspace*{1ex}\\ 
The full spectrum of the algebra of strongly regular functions has been used by M. Kashiwara in \cite{Kas} for his  
infinite dimensional algebraic geometric approach to the flag variety of a Kac-Moody group. In contrary to Kac and Peterson, he constructs the 
algebra of strongly regular functions without using the minimal Kac-Moody group. He uses the coalgebra structure of the universal enveloping 
algebra of $\g$, to construct the algebra of strongly regular functions as a certain subalgebra of the corresponding dual algebra. 
Kashiwara defines an open subscheme of the full spectrum, which has a countable covering by suitably defined big cells. To obtain his flag variety, 
Kashiwara factors the subscheme by the action of $B_f\times \{1\}$.\vspace*{1ex}\\
Taking $\F=\C$, there exists a real unitary form $K$ of $G_m$, which coincides with the compact form in the classical case. 
D. Pickrell conjectures in \cite{Pic} certain $K$-biinvariant resp. $K$-invariant measures for Kac-Moody groups. He proves the 
existence of these measures in the affine case.\\ Important for the construction of the biinvariant measures is a 
$G_f\times G_f$-space $G_f\times_{G_m}G_f$. He equips this space with a proalgebraic complex manifold structure, using a covering of big cells as an 
atlas. He also equips $G_f\times_{G_m}G_f$ with an algebra of matrix coefficients, isomorphic to the algebra of strongly regular functions. 
Important for the construction of the invariant measures is the flag $G_f$-space $\{1\}\times B_f\,\backslash\, G_f\times_{G_m}G_f$, 
also equipped with a proalgebraic manifold structure.\vspace*{1ex}\\
It is not difficult to see, that the $\C$-valued points of the subscheme of Kashiwara, as well as the $G_f\times G_f$-space $G_f\times_{G_m}G_f$ 
of Pickrell identify with the biggest $G_f\times G_f$-orbit of $\Spm\CK{G_m}$.
In a subsequent paper, we will investigate the completed flag varieties in the setting of Kashiwara and in the setting of Pickrell. Presumably the 
extended Bruhat order of \cite{M2} will be important.\vspace*{1ex}\\
Pickrell also proved a Birkhoff decomposition for $G_f\times_{G_m}G_f$. He noted, that in the affine case, many interesting
completions of $G_m$ corresponding to loop groups are embedded in the set $G_f\times_{G_m}G_f$. The Birkhoff decomposition of 
$G_f\times_{G_m}G_f$ induces the Birkhoff decomposition of these completions.\\
Now certain functional analytical closures, for example the closure used by Peterson in \cite{KP4} for his KAK-decomposition, are 
also embedded in the full set $\Spm\CK{G_m}$. Hopefully this will help to make their structure more explicit.
%
%
%
\tableofcontents\mb{}\\
%
%
%
\section{Preliminaries}
In this section we collect some basic facts about Kac-Moody algebras, minimal and formal Kac-Moody groups, the algebra of
strongly regular functions, and the corresponding monoid completion, which are used later.\\
One aim is to introduce our notation. Another aim is to put these things, which can be found in the literature, on equal footing 
appropriate for our goals.\\
The minimal Kac-Moody group, given in \cite{KP1}, \cite{KP3}, corresponds to the derived Kac-Moody algebra. We work with a 
slightly enlarged group corresponding to the full Kac-Moody algebra as in \cite{Ti1}, \cite{MoPi}. The algebra of 
strongly regular functions on this group is slightly larger, than the algebra of \cite{KP2}. We introduce this 
algebra as the restricted coordinate ring of a monoid. The formal 
Kac-Moody group has been constructed in \cite{Sl1}, starting with a realization, glueing parabolic subgroups of finite 
type, which are equipped with a proalgebraic structure. We only need the formal Kac-Moody group corresponding to a simply connected
minimal free realization, and we introduce this group by a representation theoretic construction.\\
All the material stated in this subsection about Kac-Moody algebras can be found in the books \cite{K} (most results also valid 
for a field of characteristic zero with the same proofs), \cite{MoPi}, about the minimal Kac-Moody group in \cite{KP1}, \cite{KP3}, 
\cite{MoPi}, about the formal Kac-Moody group in \cite{Sl1}, about the algebra of strongly regular functions in \cite{KP2}, 
about the faces of the Tits cone in \cite{Loo}, \cite{Sl1}, \cite{M},
and about the monoid completion of the minimal Kac-Moody group in \cite{M}.\vspace*{2ex}\\ 
We denote by $\N=\Z^+$, $\Q^+$, resp. $\R^+$ the sets of strictly positive numbers of $\Z$, $\Q$, resp. $\R\,$,
and the sets $\N_0=\Z^+_0$, $\Q^+_0$, $\R^+_0$ contain, in addition, the zero.\vspace*{0.5ex}\\
In the whole paper, $\F$ is a field of characteristic 0 and $\F^\times$ its group of units.\vspace*{1ex}\\
%
%
%
%
%
{\bf Generalized Cartan matrices:}
%
%
%
Starting point for the construction of a Kac-Moody algebra, and its associated simply connected minimal and formal Kac-Moody groups is a 
{\it generalized Cartan matrix}, which is a matrix $A=(a_{ij})\in M_{n}(\Z)$ with $a_{ii}=2$, $a_{ij}\leq 0$ for all $i\neq j$, and $a_{ij}=0$ 
if and only if $a_{ji}=0$. Denote by $l$ the rank of $A$, and set $I:=\{1,2,\ldots, n\}$.\\ 
For the properties of the generalized Cartan matrices, in particular their classification, we refer to the book \cite{K}. In this paper 
we assume $A$ to be symmetrizable.\vspace*{1ex}\\
%
%
%
{\bf Realizations:}
%
%
%
A {\it simply connected minimal free realization} of $A$ consists of dual free  $\mathbb{Z}$-modules $H$, $P$ of rank $2n-l$, 
and linear independent sets $\Pi^\vee =\{h_1,\ldots, h_n\}\subseteq H $, $\Pi=\{\al_1,\ldots,\al_n\}\subseteq P$ such that 
$\al_i(h_j)=a_{ji}\,$, $i,j=1,\dots, n$. Furthermore there exist (non-uniquely determined) fundamental do\-mi\-nant weights
$\La_1,\ldots, \La_n\in P$ such that $\La_i(h_j) =\de_{ij}$, $i,j=1,\ldots, n$.\\
$P$ is called the {\it weight lattice}, and $Q:=\Z\mb{-span}\Mklz{\al_i\,}{\,i\in I}$ the {\it root lattice}.\\
Set $Q^\pm_0:=\Z^\pm_0\mb{-span}\Mklz{\al_i\,}{\,i\in I}$, and $Q^\pm:=Q^\pm_0\setminus\{0\}$.\\
We fix a system of fundamental dominant weights $\La_1,\ldots, \La_n$, and extend $h_1,\ldots, h_n\in H$, 
$\La_1,\ldots, \La_n\in P$ to a pair of dual bases $h_1,\ldots, h_{2n-l}\in H$, $\La_1,\ldots, \La_{2n-l}\in P$. We
set $H_{rest}:=\Z \mb{-}span\Mklz{h_i}{i=n+1,\ldots,2n-l}$.\vspace*{1ex}\\
%
%
%
%
%
%
%
{\bf The Weyl group, the Tits cone and its faces:}
%
%
%
%
Identify $H$ and $P$ with the corresponding sublattices of the following vector spaces over $\F\,$:
\begin{eqnarray*}
   \h  \;\,:=\;\,  \h_\F    \;\,:=\;\,   H \otimes_{\mathbb Z} \F   &\quad,\quad &
   \h^* \;\,:=\;\, \h^*_\F  \;\,:=\;\,  P \otimes_{\mathbb Z} \F\;\;.
\end{eqnarray*}
$\h^*$ is interpreted as the dual of $\h$. Order the elements of $\h^*$ by $\la\leq\la'$ if and only if $\la'-\la\in Q_0^+$.\\ 
Choose a symmetric matrix $B\in M_n (\Q)$ and a diagonal matrix 
$D=\mb{diag}(\eps_1,\ldots, \eps_n)\,$, $\eps_1,\ldots,\eps_n\in\Qp$, such that $A=DB$.
Define a nondegenerate symmetric bilinear form on $\h$ by:
\begin{eqnarray*}
 \iB{h_i}{h} \;=\; \iB{h}{h_i} \;:=\; \al_i(h)\,\eps_i  &\qquad & i\in I\,,
 \quad h\in\h\;,\\
 \iB{h'}{h''} \;:=\; 0 \qquad\qquad &\qquad & h',\,h'' \in \h_{rest}:=H_{rest}\otimes\F \;.
\end{eqnarray*}
Denote the induced nondegenerate symmetric form on $\h^*$ also by $\iB{\;}{\;}$.\vspace*{1ex}\\
The {\it Weyl group} $\We=\We(A)$ is the Coxeter group with generators $\sigma_i\,$, $i\in I$, and
relations
\begin{eqnarray*}
       \sigma_i^2 \;=\; 1 \qquad (i\in I)    \;\:&\;,\;&\;\:
       {(\sigma_i\sigma_j)}^{m_{ij}} \;=\; 1 \qquad (i,j\in I,\,i\ne j)\;\;.
\end{eqnarray*}
The $m_{ij}$ are given by:
    $\quad \begin{tabular}{c|ccccc}
      $a_{ij}a_{ji}$ & 0 & 1 & 2 & 3 &  $\geq$ 4  \\[0.5ex] \hline 
         $m_{ij}$    & 2 & 3 & 4 & 6 &  no relation between $\sigma_i$ and $\sigma_j$ 
     \end{tabular}$\vspace*{1ex}\\
The Weyl group $\We$ acts faithfully and contragrediently on $\h$ and $\h^*$ by 
\begin{eqnarray*}
  \sigma_i h \;:=\; h - \al_i\,(h) h_i & \qquad i\in I,\quad h\in \h \;\;,\\
  \sigma_i \la \;:=\; \la - \la(h_i)\,\al_i & \qquad i\in I,\quad \la\in \h^*\;\;,
\end{eqnarray*}
leaving the lattices $H$, $Q$, $P$, and the forms invariant.\\
$\Delta_{re}:=\We\Mklz{\al_i}{i\in I}$ is called the set of {\it real roots}, and 
$\Delta_{re}^\vee:=\We\Mklz{h_i}{i\in I}$ the set of {\it real 
coroots}. The map $\al_i\mapsto h_i\,$, $i\in I$, can be extended to a $\We$-equivariant bijection
$\al \mapsto  h_\al$.\vspace*{1ex}\\
To illustrate the action of $\We$ on $\h^*_\R$ geometrically, for $J\subseteq I$ define 
\begin{eqnarray*}
  F_J &:=& \Mklz{\la\in\h^*_\R}{\la(h_i)\,=\,0 \;\mb{ for }\; i\in J\,,\;\;\;
  \la(h_i)\,>\,0\; \mb{ for }\;i\in I\setminus J}\;\;,\\
  \overline{F_J} &:=& \Mklz{\la\in\h^*_\R}{\la(h_i)\,=\,0\; \mb{ for }\;i\in J
  \,,\;\;\;\la(h_i)\,\geq\,0\; \mb{ for }\;i\in I\setminus J}\;\;. 
\end{eqnarray*}
$\overline{F_J}$ is a finitely generated convex cone with relative interior $F_J$. The parabolic subgroup $\We_J$ of $\We$ is the stabilizer of every 
element $\la\in F_J$. For $\sigma\in\We$ call $\sigma F_J$ a {\it facet} of {\it type} $J$.\\   
The {\it fundamental chamber} $\overline{C} \,:=\, \Mklz{\la\in\h^*_\R}{\la(h_i)\,\geq\,0\; \mb{ for }\;i\in I}$ is a fundamental 
region for the action of $\We$ on the convex cone $X:=\We\,\overline{C}\,$, which is called the {\it Tits cone}. The partition
$\overline{C}=\dot{\bigcup}_{J\subseteq I} F_J$ induces a $\We$-invariant partition of $X$ into facets.\\ 
A set $\Th\subseteq I$ is called {\it special}, if either $\Th=\emptyset$, or else all connected components of the 
generalized Cartan submatrix $(a_{ij})_{i,j\in\Th}$ are of non-finite type. Set 
$\Th^\bot:=\Mklz{i\in I}{ a_{ij}=0 \mb{ for all } j\in\Th }$. Every face of the Tits cone $X$ is $\We$-conjugate to 
exactly one of the faces
\begin{eqnarray*}
  R(\Th) \;\,:=\;\, X\,\cap\,\Mklz{\la\in \h_{\R}^*}{\la(h_i)=0 \mb{ for all } i\in\Th} 
          \;\,=\;\,\We_{\Th^\bot} \overline{F}_\Th&\;,\; & \Th\;\mb{ special}\;\,.
\end{eqnarray*}
The parabolic subgroup $\We_\Th$ is the pointwise stabilizer of $R(\Th)$, and the parabolic subgroup 
$\We_{\Th\cup\Th^\bot}$ is the stabilizer of the set $R(\Th)$ as a whole. \\
The relative interior of $R(\Th)$
is given by the union of the facets 
$\sigma F_{\Th\cup\Th_f}$, where $\sigma\in\We_{\Th^\bot}$, and 
$\Th_f$ is a subset of $\Th^\bot$, which is either empty, or else for which 
all connected components of $(a_{ij})_{i,j\in\Th_f}$ are of finite type.
\vspace*{1ex}\\
%
%
%
%
{\bf The Kac-Moody algebra:} 
%
%
%
%
The {\it Kac-Moody algebra} $\g=\g(A)$ is the Lie algebra over $\F$ generated by the abelian Lie algebra $\h$ and
$2n$ elements $e_i,f_i$, ($i\in I$), with the following relations, which 
hold for any $i,j \in I$, $h \in \h\,$: 
  \begin{eqnarray*}
    \left[ e_i,f_j \right] \,=\,   \delta_{ij} h_i   \;\;,\;\;   
    \left[ h,e_i \right]   \,=\,   \al_i(h) e_i    \;\;,\;\;  
    \left[ h,f_i \right]   \,=\,  -\al_i(h) f_i \;\;,    \\
    \left(ad\,e_i\right)^{1-a_{ij}}e_j  \,=\, \left(ad\,f_i\right)^{1-a_{ij}}f_j  \,=0  \qquad (i\neq j)\;\;.
  \end{eqnarray*}
The {\it Chevalley involution} $*:\g\to\g$ is the involutive anti-automorphism determined by $ e_i^*=f_i$, $ f_i^*=e_i$, $h^*=h$, 
($i\in I$, $h\in \h$).\vspace*{0.5ex}\\     
The nondegenerate symmetric bilinear form \mb{( $|$ )} on $\h$ extends uniquely to a nondegenerate symmetric invariant 
bilinear form \mb{( $|$ )} on $\g$.
We have the {\it root space decomposition}
\begin{eqnarray*}
  \g=\bigoplus_{\al \in \h^*}\g_{\al} \quad \mb{where} \quad 
  \g_\al := \Mklz{x\in \g}{[h,x]=\al(h)\,x\;\mb{ for all }\;h\in \h}\;\;.           
\end{eqnarray*}
In particular $\g_0 = \h$, $\g_{\al_i}=\F e_i$, and $\g_{-\al_i}=\F f_i$, $i\in I $.\\
The set of roots $\W:=\Mklz{\al\in\h^*\setminus\{0\}}{\g_\al\ne \{0\}}$ is 
invariant under the Weyl group, $\W=-\W$, and $\W$ spans the root lattice $Q\,$. We have $\rW\subseteq \W\,$, and 
$\iW:=\W\setminus\rW$ is called the set of {\it imaginary roots}.\\
$\W$, $\rW$, and $\iW$ decompose into the disjoint union of the sets of 
{\it positive} and {\it negative} roots $\W^\pm:=\W\cap Q^\pm$, $\rW^\pm:=\rW\cap Q^\pm$, $\iW^\pm:=\iW\cap Q^\pm$.\\ 
There is the {\it triangular decomposition} $\g = \nm \oplus \h \oplus \np $, where $\n^\pm:=\bigoplus_{\al\in {\Delta}^\pm} \g_\al $.\vspace*{1ex}\\        
%
%
%
%
%
{\bf Irreducible highest weight representations:}
%
%
%
%
For every $\La\in\h^*$ there exists, unique up to isomorphism, an irreducible representation $(L(\La),\pi_\La)$
of $\g$ with highest weight $\La$. It is $\h$-dia\-gon\-al\-iz\-able, and we denote its set 
of weights by $P(\La)$.\\
Any such representation carries a nondegenerate symmetric bilinear form $\kBl:L(\La)\times L(\La)\to\F$
which is contravariant, i.e., $\kB{v}{xw}=\kB{x^*v}{w}$ for all $v,w\in L(\La)$, $x\in\g$. This form is unique 
up to a nonzero multiplicative scalar.\vspace*{1ex}\\  
%
%
%
%
%
%
{\bf The minimal and the formal Kac-Moody group:}
%
%
%
%
We say that a Lie algebra ${\bf l}$ acts {\it locally nilpotent} on an ${\bf l}$-module $V$, if for every $v\in V$, there exists 
a positive integer $m\in\N$, such that for all $x_1,\,x_2,\ldots,x_m\in {\bf l}$ we have $x_1x_2\cdots x_m v=0$.\vspace*{0.5ex}\\
Call a $\g$-module $V$ {\it m-admissible}, if $V$ is $\h$-diagonalizable with set of weights $P(V)\subseteq P$, and $\g_\al$ acts 
locally nilpotent on $V$ for all $\al\in\rW$.\\
Examples are the adjoint representation ($\g\,$, $ad\,$), and the irreducible highest weight representations 
($L(\La)$, $\pi_\La$), $\La\in P^+:=P\cap\Cq$.\vspace*{0.5ex}\\
(Note that m-admissible is slightly different from {\it integrable}, which means $V$ is $\h$-diagonalizable, and $\g_\al$ acts 
locally nilpotent on $V$ for all $\al\in\rW$. The weights of an integrable module can be contained in 
$\Mklz{\la\in\h^*}{\la(h_i)\in\Z,\;i=1,\,\ldots,\,n}$. If the generalized Cartan matrix is degenerate, then this set is no lattice.)\vspace*{0.5ex}\\
The {\it minimal Kac-Moody group} $G=G_m=G_m(A)$ can be characterized in the following way:\\
$\bullet$ The group $G$ acts on every m-admissible representation. Two elements $g,g'\in G$ are equal if and only if  
for all m-admissible modules $V$, and for all $v\in V$, we have $gv=g'v$.\\
$\bullet$ (1) For every $h\in H$, $s\in\F^\times$ there exists an element $t_h(s)\in G$, such that for any m-admissible 
representation $(V,\pi)$ we have
\begin{eqnarray*}
  t_h(s)v_\la &=& s^{\la(h)}v_\la \quad,\quad v_\la\in V_\la\;\,,\;\,\la\in P(V)\;.
\end{eqnarray*}
(2) For every $x\in\g_\al$, $\al\in\rW$, there exists an element $\exp(x)\in G$, such that for 
any m-admissible representation $(V,\pi)$ we have
\begin{eqnarray*}
  \exp(x)v &=& \exp(\pi(x)) v \quad,\quad v\in V\;.
\end{eqnarray*}
$G$ is generated by the elements of (1) and (2).\vspace*{1ex}\\
Call a $\g$-module $V$ {\it f-admissible}, if $V$ is m-admissible, and $\n^+$ acts locally nilpotent on $V$.\\
Examples are the representations ($L(\La)$, $\pi_\La$), $\La\in P^+=P\cap\Cq$.\vspace*{0.5ex}\\
Set $\n_f:=\prod_{\al\in \Delta^+}\g_\al$ and $\g_f:=\n^-\oplus\h\oplus\n_f$. The Lie bracket of $\g$ extends in the obvious
way to a Lie bracket of $\g_f$. Every f-admissible $\g$-module can be extended to a $\g_f$-module. The Lie algebra $\g_f$ should be 
interpreted as the Lie algebra of the {\it formal Kac-Moody group} $G_f=G_f(A)$, which can be characterized in the following way:\\
$\bullet$ The group $G_f$ acts on every f-admissible representation. Two elements $g,g'\in G_f$ are equal if and only if  
for all f-admissible modules $V$, and for all $v\in V$ we have $gv=g'v$.\\  
$\bullet$ (3) $G_f$ contains $G$.\\
(4) For every $x\in\n_f$ there exists an element $\exp(x)\in G_f$, such that for any f-admissible representation $(V,\pi)$ we have
\begin{eqnarray*}
  \exp(x)v &=& \exp(\pi(x)) v \quad,\quad v\in V\;\;.
\end{eqnarray*}
$G_f$ is generated by $G$ and the elements of (4).\vspace*{0.5ex}\\
The $\g$-module $\g_f$ is not f-admissible.  Nevertheless $G_f$ acts on $\g_f$, extending the adjoint action of $G$ on $\g$, 
compare \cite{Sl1}, Section 5.11.\vspace*{1ex}\\
Both Kac-Moody groups act faithfully on $\bigoplus_{\La\in P^+}L(\La)$. They have the following important structural 
properties:\vspace*{0.5ex}\\
1) The elements of (1) induce an embedding of the torus $H\otimes_\Z\F^\times$ into $G\subseteq G_f$. Its image is denoted 
by $T$.\\
For $\al\in\rW$ the elements of (2) induce an embedding of $(\g_\al,+)$ into $G\subseteq G_f$. Its image $U_\al$ is called the 
root group belonging to $\al$.\\
Let $\al\in \prW$ and $x_\al\in\g_{\al}$, $x_{-\al}\in\g_{-\al}$ such that $[x_\al,x_{-\al}]=h_\al$. There exists an 
injective homomorphism of groups $\phi_\al:\,\mb{SL}(2,\F) \to G$ with 
\begin{eqnarray*}
 \phi_\al\left(\begin{array}{cc}
  1 & s\\
  0 & 1
  \end{array}\right) \;:=\; \exp(s x_\al)  \;\,,\,\;
   \phi_\al\left(\begin{array}{cc}
  1 & 0\\
  s & 1
  \end{array}\right) \;:=\; \exp(sx_{-\al})\;\,,\,\;(s\in \F^\times)\;.
\end{eqnarray*}
2) Denote by $N$ the subgroup generated by $T$ and $
  n_\al := \phi_\al\left(\begin{array}{cc}
        0      & 1  \\
  -1 & 0
  \end{array}\right)$, $\al\in \Delta_{re}$. The Weyl group $\We$ can be identified with the group $N/T$ by the isomorphism 
$\kappa:\,\We \to  N/T$, which is given by $\kappa(\sigma_\al):=n_\al(1) T$, $\al\in \rW $.\\
We denote an arbitrary element $n\in N$ with $\kappa^{-1}(nT)=\sigma\in\We$  by $n_\sigma$. The set of weights $P(V)$ of an
m-admissible $\g$-module $(V,\pi)$ is $\We$-invariant, and $ n_\sigma V_\la =V_{\sigma \la}$, $\la\in P(V)$.
\vspace*{1ex}\\
3) Let $U^\pm$ be the subgroups generated by $U_\al$, $\al\in\Delta_{re}^\pm $. Let $U_f:=\exp(\n_f)$. Then $U^\pm$ and $U_f$ are
normalized by $T$. Set
\begin{eqnarray*}
   B^\pm\;\,:=\,\; T\ltimes U^\pm &\quad,\quad &  B_f\;\,:=\,\; T\ltimes U_f\;\;.
\end{eqnarray*}
The pairs ($B^\pm$, $N$) are  twinned BN-pairs of $G$ with the property $B^+\cap B^-= B^\pm\cap N = T$. The pair 
($B_f$, $N$) is a BN-pair of $G_f$ with $B_f\cap N=T$. We have the following decompositions, called  {\it Bruhat} and 
{\it Birkhoff decompositions}: 
\begin{eqnarray*}
  G \;\,=\;\, \dot{ \bigcup_{\sigma\in{\cal W}}} B^\epsilon\sigma B^\delta \quad,\quad 
  G_f \;\,=\;\, \dot{ \bigcup_{\sigma\in{\cal W}}} B^\epsilon\sigma B_f \qquad,
  \qquad \epsilon,\delta\;\in\;\{\,+\,,\,-\,\}\;\;.
\end{eqnarray*}
4) There are also Levi decompositions of the standard parabolic subgroups. In this paper we only use the corresponding
decompositions for the groups $U^\pm$ and $U_f$:
Set $\W^\pm_J:=\W^\pm\cap \sum_{j\in J}\Z\,\al_j$, and ${(\W^J)}^\pm:=\W^\pm\setminus\sum_{j\in J}\Z\,\al_j$. Similarly define 
$(\W_J)^\pm_{re}$ and $(\W^J)^\pm_{re}$ by replacing $\W^\pm$ by $\W_{re}^\pm$. Set 
$(\n_J)^\pm:=\bigoplus_{\al\in \Delta_J^\pm}\g_\al$, $(\n_f)_J:=\prod_{\al\in\Delta_J^+}\g_\al$, and 
$(\n_f)^J:=\prod_{\al\in(\Delta^J)^+}\g_\al$. We have:
\begin{eqnarray*}
   U^\pm\;\,=\;\, U^\pm_J \ltimes (U^J)^\pm  &\quad,\quad &  U_f\;\,=\;\, (U_f)_J \ltimes (U_f)^J\;\;.
\end{eqnarray*}
Here $U_J^\pm$ is the group generated by the root groups $U_\al$, $\al\in (\W_J)^\pm_{re}$. $(U^J)^\pm$ is the smallest normal subgroup of $U^\pm$, 
containing the root groups $U_\al$, $\al\in (\W^J)^\pm_{re}$. This group equals $\,\bigcap_{\sigma\in{\cal W}_J}\sigma U^\pm\sigma^{-1}$.
Furthermore $(U_f)_J:=\exp((\n_f)_J)$ and $(U_f)^J:=\exp((\n_f)^J)$.\vspace*{1ex}\\ 
The derived minimal Kac-Moody group $G'$ is identical with the Kac-Moody group as defined in \cite{KP1}. It is generated by 
the root groups $U_\al$, $\al\in \rW$. We have $G = G'\rtimes T_{rest}$, where $T_{rest}:=H_{rest}\otimes_{\mathbb Z} \F$ is a 
subtorus of $T$.\\
The group $G_f$ is identical with the Kac-Moody group of \cite{Sl1} for a simply connected minimal free realization.
\vspace*{1ex}\\
%
%
%
{\bf The monoid $\GD$:}
%
%
%
The category $\cal O$ is defined as follows: Its objects are the $\g$-modules $V$, which have the properties:\\
(1) $V$ is $\h$-diagonalizable with finite dimensional weight spaces.\\
(2) There exist finitely many elements $\la_1,\,\ldots,\,\la_m\in\h^*$, such that the set of weights $P(V)$ of $V$ is contained in the union 
$\bigcup_{1=1}^m D(\la_i)$, where $D(\la_i):=\Mklz{\la\in\h^*}{\la\leq \la_i}$.\\
The morphisms of $\cal O$ are the morphisms of $\g$-modules.\vspace*{1ex}\\ 
For $V$ a module of ${\cal O}$ and $v\in V$, we denote by $supp(v)$ the set of weights of the nonzero weight space componsnts of $v$.\vspace*{1ex}\\
A module of ${\cal O}$ is m-admissible if and only if it is f-admissible, and we call such a module {\it admissible}. We denote by 
${\cal O}_{adm}$ the full subcategory of the category $\cal O$, whose objects are admissible modules.\\ 
There is a complete reducibility theorem. Every object of ${\cal O}_{adm}$ is isomorphic to a direct sum of the admissible irreducible highest weight 
modules $L(\La)$, $\La\in P^+$.\vspace*{1ex}\\
The set of weights of a module of ${\cal O}_{adm}$ is contained in $X\cap P$, because we have
\begin{eqnarray*}
  \bigcup_{\La\in P^+} P(\La) &=& X\cap P\;\;.
\end{eqnarray*}
Let $\La\in P^+$, and $\Th$ be special. Because the set of weights $P(\La)$ is contained in the convex hull of $\We\La\subseteq \h_\R^*$, 
we find easily 
\begin{eqnarray*}
   P(\La)\cap R(\Th) \;\,=\;\,\emptyset &\mb{ if and only if }& \La \notin R(\Th)\;\;.
\end{eqnarray*}
The monoid $\GD$ can be characterized in the following way:\\
$\bullet$ The monoid $\GD$ acts on every module of ${\cal O}_{adm}$. Two elements $\hat{g},\hat{g}'\in \GD$ are equal if and only if  
for all modules $V$ of ${\cal O}_{adm}$, and for all $v\in V$, we have $\hat{g}v=\hat{g}'v$.\\ 
$\bullet$
(1) $\GD$ is an extension of the minimal Kac-Moody group $G$.\vspace*{0.5ex}\\
(2) For every face $R$ of the Tits cone there exists an element $e(R)\in \GD$, such that for every module $V$ of ${\cal O}_{adm}$ we have
\begin{eqnarray*}
  e(R)v_\la &=& \left\{ \begin{array}{ccc}
    v_\la &\;&\la\in R\\
      0   &\;& \mb{else}
  \end{array}\right. \quad,\quad \;v_\la\in V_\la\,,\; \;\la\in P(V)\;\;.
\end{eqnarray*}
$\GD$ is generated by $G$ and the elements of (2).\vspace*{1ex}\\
Note that the monoid $\GD$ acts faithfully on the sum $\bigoplus_{\La\in P^+}L(\La)$.\vspace*{1ex}\\  
The {\it Chevalley involution} $*:\,\GD \to \GD$ is the involutive anti-isomorphism determined by $\exp(x_\al)^*:=\exp(x_\al^*)$, 
$t^*:= t$, $e(R)^* := e(R)$, where $x_\al\in\g_\al$, $\al\in\rW$, $t\in T$, and $R\in\RkX$.\\
It is compatible with any nondegenerate symmetric contravariant form $\kBl$ on any module $V$ of ${\cal O}_{adm}$, i.e., $\kB{xv}{w}=\kB{v}{x^*w}$, 
$v,w\in V$, $x\in\GD$.\vspace*{1ex}\\
The following formulas are useful for computations in $\GD$:\vspace*{0.5ex}\\
$\bullet$ Let $R$, $S$ be faces of the Tits cone, and $n_\sigma\in N$. Then
\begin{eqnarray*}
 e(R)e(S)=e(R\cap S) \qquad ,\qquad  n_\sigma e(R) n_\sigma^{-1}\;\,=\;\, e(\sigma R)\;\;.
\end{eqnarray*}  
$\bullet$ An element $g$ of $T$, $N$, $U$, $U^-$, resp. $G$ satisfies 
\begin{eqnarray*}
   e(R(\Th)) g &=& e(R(\Th))
\end{eqnarray*}
if and only if it satisfies
\begin{eqnarray*}
    g^* e(R(\Th))  &=& e(R(\Th))
\end{eqnarray*}
if and only if it is contained in $T_\Th$, $N_\Th$, $U_\Th$, $U_\Th^-\ltimes (U^{\Th\cup\Th^\bot})^{-}$, resp. $G_\Th\ltimes U^{\Th\cup\Th^\bot}$.
Here $T_\Th$ is the subtorus of $T$ generated by $t_{h_j}(s)$, $j\in \Th$, $s\in\F^\times$, $N_\Th$ is the subgroup of $N$ generated by $T_\Th$ and 
$n_{\al_j}$, $j\in \Th$, and $G_\Th$ is the subgroup of $G$ generated by $U_{\al_j}^\pm$, $j\in \Th$.\vspace*{1ex}\\
$\bullet$ An element $g$ of $T$, $N$, $U$, $U^-$, resp. $G$ satisfies 
\begin{eqnarray*}
   g e(R(\Th)) g^{-1} &=& e(R(\Th))
\end{eqnarray*}
if and only if it is contained in the groups $T$, $ N_{\Th\cup\Th^\bot} T$,  $U_{\Th\cup\Th^\bot}$, $U_{\Th\cup\Th^\bot}^-$, resp. 
$G_{\Th\cup\Th^\bot} T$.\vspace*{1ex}\\ 
$\bullet$ In particular we have
\begin{eqnarray*}
    U e(R(\Th))  &=& U_{\Th^\bot} e(R(\Th))\;\,=\;\, 
    e(R(\Th))U_{\Th^\bot}\;\;,\\
     e(R(\Th))U^-  &=& e(R(\Th))U^-_{\Th^\bot}\;\,=\;\, 
      U^-_{\Th^\bot}e(R(\Th))\;\;.
\end{eqnarray*}
The minimal Kac-Moody group $G$ is the unit group of $\GD$. Every idempotent is $G$-conjugate to some idempotent $e(R(\Th))$, $\Th$ special. We have
\begin{eqnarray*}
  \GD&=& \dot{\bigcup_{\Th\,special}} G e(R(\Th)) G\;\;.
\end{eqnarray*}
The Weyl group acts on the monoid (\,$\RkX\,,\,\cap$\,). The semidirect product $\RkX\rtimes\We$ consists of the set $\RkX\times\We$ with the 
structure of a monoid given by 
\begin{eqnarray*}
  (R,\sigma)\cdot(S,\tau) &:=& (R\cap\sigma S,\sigma \tau)\;\;.
\end{eqnarray*}
For $R\in\RkX$ let $Z_{\cal W}(R):=\Mklz{\sigma\in \We}{\sigma\la=\la \mb{ for all }\la\in R}$ be the pointwise stabilizer of $R$. The Weyl monoid 
$\WeD$ is defined as the monoid $\RkX\rtimes\We$ factored by the congruence relation
\begin{eqnarray*}
   (R,\sigma) \sim (R',\sigma') &:\iff & R\:=\:R'\quad\mb{and}\quad\sigma'\sigma^{-1}\in Z_{\cal W}(R)\;\;.
\end{eqnarray*}
We denote the congruence class of $(R,\sigma)$ by $\ve{R}\sigma$.\\
Assigning to $\sigma\in\We$ the element $\sigma:=\sigma e(X)\in\WeD$, the Weyl group $\We$ identifies with the unit group of $\WeD$. The partition of 
$\WeD$ into $\We\times\We$-orbits is given by
\begin{eqnarray*}
  \WeD &=& \dot{\bigcup_{\Th \;special}} \,\We\,\ve{R(\Th)}\,\We\;\;.
\end{eqnarray*}
Assigning to $e(R)\in \RkX$ the element $\ve{R}:=\ve{R}1\in\WeD$, the monoid $(\,\RkX\,,\,\cap \,)$ embedds into $\WeD$. Its image are the idempotents of 
$\WeD$.\\ 
For $J\subseteq I$ denote by $\We^J$ the minimal coset representatives of $\We/\We_J$, and denote by $\mb{}^J\We$ the minimal coset 
representatives of $\We_J\backslash \We$. It is easy to see that there are the following uniquely determined normal forms of an element 
$\hat{\sigma}\in \WeD$:
\begin{eqnarray*}
   \hat{\sigma}\;\,=\;\, \sigma_1 \ve{R(\Th)} \sigma_2 &\mb{ with }&
   \Th \mb{ special }\,,\;\sigma_1\in \We^{\Th\cup\Th^\bot}\,,\;
   \sigma_2\in\mb{}^{\Th}\We\;.  \\
   \hat{\sigma}\;\,=\;\, \tau_1 \ve{R(\Th)} \tau_2 &\mb{ with }&
   \Th \mb{ special }\,,\;\tau_1\in \We^{\Th}\,,\;
   \tau_2\in\mb{}^{\Th\cup\Th^\bot}\We \;.
\end{eqnarray*}
Let $J\subseteq I$. We call the submonoid $\WeD_J$, which is generated by $\We_J$ and the elements $\ve{R(\Th)}$, $\Th\subseteq J$ special, a 
parabolic submonoid. Denote by $J^\infty$ the union of all connected components of nonfinite type of $J$. We have
\begin{eqnarray*}
  \WeD_J \;\,=\;\, \dot{\bigcup_{R\,a\,face\,of\,X\atop R\supseteq R(J^\infty)}} \We_J \,\ve{R}\;\,=\;\, 
  \dot{\bigcup_{\Xi \,special\atop  \,\Xi\subseteq J^\infty}} \We_J \,\ve{R(\Xi)}\,\We_J\;\;.
\end{eqnarray*}
We get an abelian submonoid of $\GD$ by $\TD:=\dot{\bigcup}_{R\in{\cal R}(X)} T e(R) $. We get a submonoid of $\GD$ by 
$\ND:=\dot{\bigcup}_{R\in{\cal R}(X)} N e(R)$. Define a congruence relation on $\ND$ as follows:
\begin{eqnarray*}
\quad\hat{n}\:\sim\:\hat{n}'\quad:\iff\quad\hat{n}T\:=\:\hat{n}'T\quad\iff \quad\hat{n}'\:\in\:\hat{n}T\quad\iff\quad\hat{n}\in\hat{n}'T
\end{eqnarray*}
The Weyl monoid $\WeD$ is isomorphic to the monoid $\ND/T$, an isomorphism $\kappa:\WeD\to \ND/T$ given by $\kappa(\sigma\ve{R})=n_\sigma e(R) T$.
\vspace*{1ex}\\
$\GD$ has {\it Bruhat} and {\it Birkhoff decompositions}: 
\begin{eqnarray*}
  \GD &=& \dot{ \bigcup_{\hat{n}\in \widehat{N}}}\,\; U^\epsilon\,\hat{n}\, U^\delta \;\,=\;\,
   \dot{ \bigcup_{\hat{\sigma}\in \widehat{\cal W}}}\,\; B^\epsilon\,\hat{\sigma}\, B^\delta \quad,\quad  \epsilon,\,\delta \;\in \;\{\,+\,,\,-\,\}\;\;.
\end{eqnarray*}
%
%
%
%
%
{\bf The coordinate ring of $\GD$, and the algebra of strongly regular functions:}
%
%
%
%
For a module $V$ of ${\cal O}_{adm}$, $v,w\in V$, and $\kBl$ a nondegenerate symmetric contravariant bilinear form on $V$, call the function 
$\ti{f}_{vw}:\,\GD\to\F$ defined by $\ti{f}_{vw}(x):=\kB{v}{xw}$, $x\in\GD$, a {\it matrix coefficient} of $\GD$. The set of all such matrix 
coefficients $\FK{\GD}$ is an algebra. It is an integral domain, and admits a {\it Peter-Weyl theorem}: Equip $\FK{\GD}$ with an action $\pi$ 
of $\GD\times \GD$, and an involutive automorphism $*$ by: 
\begin{eqnarray*}
  \begin{array}{ccc}
    (\pi(g,h)\,f)(x) &:=& f(g^* x h) \\
     f^* (x) &:=& f(x^*)
   \end{array} &, & g,x,h\in \GD\,,\quad f\in\FK{\GD}\;\;.
\end{eqnarray*} 
For every $\La\in P^+$ fix a nondegenerate symmetric contravariant bilinear form on $L(\La)$. The map 
$\bigoplus_{\La\in P^+} L(\La)\otimes L(\La) \to  \FK{\GD}$ induced by $v\otimes w\mapsto \ti{f}_{vw}\,$ 
is an isomorphism of $\GD\times \GD$-modules. It identifies the direct sum of the switch maps of the factors with the involution.\vspace*{0.5ex}\\
The monoids $\TD$, $\ND$, $\GD$ are the Zariski closures of $T$, $N$, $G$, and $G$ is the Zariski open dense unit group of $\GD$.\vspace*{0.5ex}\\
The algebra of {\it strongly regular functions} $\FK{G}$ is obtained by restricting the functions of $\FK{\GD}$ onto $G$. 
The restriction map is an isomorphism from $\FK{\GD}$ to $\FK{G}$.\vspace*{0.5ex}\\
Restricting the functions of $\FK{G}$ onto $G'$, resp. $T_{rest}$ gives the algebras $\FK{G'}$, resp. $\FK{T_{rest}}$, the first
identical with the algebra of strongly regular functions as defined in \cite{KP2}, the second the classical 
coordinate ring of the torus $T_{rest}\,$. $\FK{G}$ is isomorphic to $\FK{G'}\otimes \FK{T_{rest}}$, by the comorphism dual to the 
multiplication map $G'\times T_{rest}\to G\,$.\vspace*{1ex}\\   
%
%
%
{\bf Substructures:}
%
%
%
For $\emptyset\neq J\subseteq I$ the submatrix $A_J:=(a_{ij})_{i,j\in J}$ of $A$ is a generalized Cartan matrix. There exist
saturated sublattices $H(A_J)\subseteq H$, $P(A_J)\subseteq P$ with $(h_j)_{j\in J}\subseteq H(A_J)$, 
$(\al_j)_{j\in J}\subseteq P(A_J)$, giving a simply connected minimal free realization of $A_J$.
We have $P=P(A_J)\oplus H(A_J)^\bot$, and the projections of $\La_j$, $j\in J$, to $P(A_J)$ are a system of fundamental dominant weights.\\
The corresponding Kac-Moody algebra $\g(A_J)$ embeds in $\g$. If we identify $\g(A_J)$ with its image, then the set of roots $\W(A_J)$ identifies 
with $\W_J:=\W\cap \sum_{j\in J}\Z\,\al_j$, and the Weyl group $\We(A_J)$ identifies with the parabolic subgroup $\We_J$.\\
The face lattice of the Tits cone of $A_J$ embeds onto a sublattice of the face lattice of the Tits cone of $A$, and the Weyl monoid 
$\widehat{\We(A_J)}$ identifies with the parabolic submonoid $\WeD_J$.\\
The minimal and formal Kac-Moody groups $G(A_J)$, $G_f(A_J)$, the monoid $\widehat{G(A_J)}$ embed in $G$, $G_f$, $\GD$ in the
obvious way.\vspace*{1ex}\\
The images of these embedding depend on the choice of the sublattice $H(A_J)$, only $H_J:=\Z\mb{-}span\Mklz{h_j}{j\in J}$ is 
uniquely determined by $A_J$. Denote by $\widehat{G'}$ the submonoid of $\GD$, which is 
generated by $G'$ and the elements $e(R)$, $R$ a face of $X$. The images of $\g(A_J)'$, $G(A_J)'$, $G_f(A_J)'$, and $\widehat{G(A_J)'}$ are
independent of this choice, and denoted by $\g_J$, $G_J$, $(G_f)_J$, and $\widehat{G}_J$.\vspace*{1ex}\\
An admissible irreducible highest weight module $L(\La)$ of $\g$, equipped with a nondegenerate contravariant symmetric bilinear form, decomposes as a 
$\g(A_J)$-module into an orthogonal direct sum of admissible irreducible highest weight modules of $\g(A_J)$, which are $\h$-invariant. In particular 
$L_J(\La):=U(\n_J^-)L(\La)_\La$ is an admissible irreducible highest weight module of $\g(A_J)$, its highest weight given by the projection of $\La$ to 
$P(A_J)$, ($\La\in P^+$).\vspace*{1ex}\\
The coordinate rings $\FK{\widehat{G(A_J)'}}$, $\FK{G(A_J)'}$ identify with the restrictions of $\FK{\GD}$, $\FK{G}$ to
$\GD_J$, $G_J$. (A similar statement for $\FK{\widehat{G(A_J)}}$, $\FK{G(A_J)}$ is not valid.)\vspace*{1ex}\\
To simplify the notation of many formulas, it is useful to set $\g_\emptyset:=\{0\}$, and $G_\emptyset:=\GD_{\emptyset}:=(G_f)_\emptyset:=\{1\}$.\vspace*{1ex}\\
For $M\subseteq \GD$ and $J\subseteq I$ set $M_J:=M\cap \GD_J$, and similarly for $M\subseteq G_f$ set $M_J:=M\cap (G_f)_J$.
%
%
%
%
%
%
\section{An easy algebraic geometric setting \label{Setting}} 
%
%
%
In this section we develop an easy algebraic geometric setting, which is useful to determine the $\F$-valued points of the algebra of strongly regular 
functions.\\
This section is a complement to \cite{M}, Section 3, but it can be read independently. The definition of a morphism in \cite{M}, Section 3, is 
more restrictive, because it has to  preserve an extra structure, which we do not need here.\vspace*{1ex}\\ 
We will consider certain nonempty sets $A$ equipped with point separating algebras of functions $\FK{A}$, which we call coordinate rings.\\
The closed sets of the Zariski topology on such a set $A$ are given by the zero sets of the functions of $\FK{A}$. Note that $A$ is irreducible 
if and only if $\FK{A}$ is an integral domain.\\
A morphism of sets with coordinate rings ($A,\FK{A}$) and ($B,\FK{B}$) consists of a map $\phi:A\to B$, whose comorphism $\phi^*:\FK{B}\to\FK{A}$ 
exists. In particular a morphism is Zariski continuous.\vspace*{1ex}\\
If ($B,\FK{B}$) is a set with coordinate ring, and $A$ is a nonempty subset of $B$, then we get a coordinate ring on $A$ by restricting the functions of 
$\FK{B}$ to $A$.\\
If ($A,\FK{A}$) is a set with coordinate ring and $f\in\FK{A}\setminus\{0\}$, the principal open set $D_A(f):=\Mklz{a\in A}{f(a)\neq 0}$ is equipped 
with a coordinate ring by identifying the localization $\FK{A}_f$ in the obvious way with an algebra of functions on $D_A(f)$. If $A$ is irreducible, then also $D_A(f)$ is irreducible.\\
If ($A,\FK{A}$) and ($B,\FK{B}$) are sets with coordinate rings, then $A\times B$ is equipped with a coordinate ring by identifying the tensor product 
$\FK{A}\otimes\FK{B}$ in the obvious way with an algebra of functions on $A \times B$. If $A$ and $B$ are irreducible, then also $A\times B$ is 
irreducible.\vspace*{1ex}\\
To construct the sets with coordinate rings, which we will use later, fix nondegenerate symmetric contravariant forms $\kBl$ on all modules 
$L(\La)$, $\La\in P^+$, and extend to a form on $\bigoplus_{\La\in P^+} L(\La)$, also denoted by $\kBl$, by requiring $L(\La)$ and $L(\La')$ to be 
orthogonal for $\La\neq\La'$.
For $v, w\in L(\La)$, $\La\in P^+$, define a linear function 
\begin{eqnarray*}
   f_{v w}:\;\Endo   &\to & \F 
\end{eqnarray*}
by $f_{vw}(\phi):=\kB{v}{\phi w}$, where $\phi\in \Endo$.\vspace*{1ex}\\
We equip the subalgebra
\begin{eqnarray*}
         \grAdj &:=& \Mklz{\,\phi\in\Endo\;}{\;\begin{array}{c} \mb{The adjoint }\;\phi^* \;\mb{ exists and}\\    
         \phi \left(L(\La)\right)\subseteq L(\La),\;\La\in P^+\;.    
         \end{array}\,}                                            
\end{eqnarray*}
of the algebra of endomorphisms with the coordinate ring $\FK{\grAdj}$, which is generated by the linear subspace
\begin{eqnarray*}
 span\Mklz{f_{v w}\res{gr\mb{-}Adj}\,}{\,v,w\in L(\La)\;\,,\,\; \La\in P^+\;}
\end{eqnarray*}
of the linear dual of $\grAdj$. This coordinate ring is isomorphic to the symmetric algebra in this subspace.\\
The multiplication map is no morphism, but the left and right multiplications with elements of $\grAdj$, and also the adjoint map $*:\grAdj\to\grAdj$ 
are morphisms.\\ 
For a set $M\subseteq \grAdj$ we denote by  $\Aa{M}$ its Zariski closure in $\grAdj$. If $M$ is a ($*$-invariant) submonoid of $\grAdj$, then also 
$\Aa{M}$ is a ($*$-invariant) submonoid. The left and right translations with elements of $\Aa{M}$ (and the map $*:\Aa{M}\to \Aa{M}$) are 
morphisms.\vspace*{1ex}\\
A function $f\in\FK{M}$ induces a function on $\Spm\FK{M}$, assigning $x\in \Spm\FK{M}$ the value $x(f)$. In this way, $\Spm\FK{M}$
is equipped with a coordinate ring isomorphic to $\FK{M}$. Its Zariski topology coincides with the relative topology induced by the
topology of the spectrum of $\FK{M}$.\\
To investigate $\Spm\FK{M}$, we introduce two new monoids $\Ae{M}$, $\Ae{M^*}$ with coordinate rings.\vspace*{1ex}\\
Equip the subalgebra 
\begin{eqnarray*}
         \grEnd &:=& \Mklz{\,\phi\in\Endo\;}{\;\phi \left(L(\La)\right)\subseteq L(\La),\;\La\in P^+\;.\,}                                            
\end{eqnarray*}
of the algebra of endomorphisms with the coordinate ring $\FK{\grEnd}$ generated by the linear subspace
\begin{eqnarray*}
 span\Mklz{f_{v w}\res{gr\mb{-}End}\,}{\,v,w\in L(\La)\;\,,\,\; \La\in P^+\;}
\end{eqnarray*}
of the linear dual of $\grEnd$, which is isomorphic to $\bigoplus_{\La\in P^+}L(\La)\otimes L(\La)$.  This coordinate ring is isomorphic to the 
symmetric algebra in this subspace. The restriction of $\FK{\grEnd}$ to $\grAdj$ coincides with the coordinate ring of $\grAdj$.\\
The right multiplication with an element $\phi\in \grEnd$ is a morphism. We can guarantie the left multiplication $l_\phi$ to be a morphism, only if we 
restrict to $\phi\in\grAdj$, because then we have
\begin{eqnarray*} 
  f_{vw}\res{gr\mb{-}End}\,\circ\: l_\phi \;\,=\;\,f_{\phi^*v\, w}\res{gr\mb{-}End}\;\;\mb{ for all }\;\; v,w\in L(\La)\;,\;\La\in P^+.
\end{eqnarray*}
For a set $M\subseteq \grEnd$ we denote by  $\Ae{M}$ its Zariski closure in $\grEnd$. If $M$ is a submonoid of $\grEnd$, then $\Ae{M}$ is not 
necessarily a submonoid of $\grEnd$. But we have:
\begin{Prop}\label{Setting1} If $M$ is a submonoid of $\grAdj$, then $\Ae{M}$ is a submonoid of $\grEnd$.
\end{Prop}
Also the following proposition is easy to prove:
\begin{Prop}\label{Setting1b} Let $M$ be a submonoid of $\grAdj$. For $\phi\in \Ae{M^*}$ and $\psi\in\Ae{M}$ we get a homomorphism of algebras 
$\pi(\phi,\psi):\FK{M}\to\FK{M}$ by
\begin{eqnarray*}
  \pi(\phi,\psi)\left(f_{v w}\res{M}\right) \;:=\; f_{\phi v\,\psi w}\res{M} &\quad,\quad & v,w\in L(\La)\;,\;\La\in P^+\;\;.
\end{eqnarray*}
$\pi$ is an action of  $\Ae{M^*}\times\Ae{M}$ on $\FK{M}$.
\end{Prop}
Due to this proposition we get an action of $\Ae{M^*}\times\Ae{M}$ on $\Spm\FK{M}$ from the right by
\begin{eqnarray*}
  \al\circ\pi(\phi,\psi)  \;\,,\quad\mb{ where }\;\,  \al\in\Spm\FK{M}\;\,,\,\;\phi\in \Ae{M^*}\;\,,\;\,\psi\in \Ae{M}\;\;.
\end{eqnarray*}
This is an action by morphisms of $\Spm\FK{M}$. But in general the map $(\Ae{M^*}\times\Ae{M})\times\Spm\FK{M}\to\Spm\FK{M}$ is no morphism. In general 
also the right translations of this map are no morphisms.\vspace*{1ex}\\
The next Proposition, which is easy to prove, describes the part of $\Spm\FK{M}$ obtained by applying $\Ae{M^*}\times\Ae{M}$ to the evaluation map 
of $\FK{M}$ in the unit of $M$.
\begin{Prop}\label{Setting2} Let $M$ be a submonoid of $\grAdj$.\\
1) For $\phi\in\Ae{M^*}$ and $\psi\in\Ae{M}$ we get a point $\al(\phi,\psi)\in\Spm\FK{M}$ by 
\begin{eqnarray*}
  \al(\phi,\psi)(f_{vw}\res{M})\;\,:=\;\,\kB{\phi v}{\psi w}  \;\,\mb{ for all }\;\; v,w\in L(\La)\,,\; \La\in P^+\;\;.
\end{eqnarray*}
If we equip $\Ae{M^*}\times\Ae{M}$ with the right action on itself, then the map $\al:\Ae{M^*}\times\Ae{M}\to\Spm\FK{M}$ is equivariant. Furthermore we 
have
\begin{eqnarray*}
   \al(\phi,z\psi)\;\,=\,\;\al(z^*\phi,\psi) &\mb{ for all }& \phi\in\Ae{M^*}\;,\;\psi\in\Ae{M}\;,\;\mb{ and }\; 
    z\in\,\Aa{M}\;\;.
\end{eqnarray*}
2) For every $\phi\in \Ae{M^*}$ the map $\al(\phi,\,\cdot\,):\,\Ae{M}\to\Spm\FK{M}$, which assigns $\psi$ the point $\al(\phi,\psi)$, is a morphism. 
The map $\al(1,\,\cdot\,)$ is injective. \\  
For every $\psi\in \Ae{M}$ the map $\al(\,\cdot\,,\psi):\,\Ae{M^*}\to\Spm\FK{M}$ which assigns $\phi$ the point $\al(\phi,\psi)$, is a morphism. The map 
$\al(\,\cdot\,,1)$ is injective.
\end{Prop}
{\bf Remarks:}\\
1) We often say $\Ae{M^*}\times\Ae{M}$ maps to $\FK{M}$, without mentioning the  map $\al$. In general this map is no morphism.\vspace*{0.5ex}\\
2) The maps $\{1\}\times \Ae{M}\to \Spm\FK{M}$, and $\Ae{M^*}\times \{1\}\to \Spm\FK{M}$ are 
injective. But if $M$ is nontrivial, then the map $\Ae{M^*}\times\Ae{M}\to \Spm\FK{M}$ is not injective:\\
Let $\sim$ be the equivalence relation on $\Ae{M^*}\times\Ae{M}$ generated by
\begin{eqnarray*}
     (z^*x,y) \,\sim\, (x,zy) &\quad,\quad & z\,\in\, \Aa{M}   \;\;.
\end{eqnarray*}
Equip the quotient set, which we denote by $\Ae{M^*}\times_{\Aa{M}}\Ae{M}$, with the induced
$\Ae{M^*}\times\Ae{M}$-action from the right. Then $\al$ factors to an equivariant map
\begin{eqnarray*}
  \ti{\al}\,:\,\Ae{M^*}\times_{\Aa{M}}\Ae{M} &\to & \Spm\FK{M}\;\;.
\end{eqnarray*}
3) The map $*:M\to M^*$ is an isomophism. Its comorphism $\mb{}^*:\FK{M^*}\to\FK{M}$ is given by
\begin{eqnarray*} 
  (f_{v w}\res{M^*})^* & = & f_{w v}\res{M}\;\,\mb{ for all } \;\,v,w\in L(\La)\;\,, \;\,\La\in P^+\;\;.
\end{eqnarray*} 
We omit to state the easy compatibility conditions with the maps $\al$, $\pi$, which correspond to $\FK{M}$, and $\FK{M^*}$.\vspace*{1ex}\\
The properties stated in the second part of the last proposition are sufficient to guarantie the irreducibility of certain orbits, belonging to the 
action of the product of a irreducible subgroup of $\Ae{M^*}$ and a irreducible subgroup of $\Ae{M}$.
\begin{Theorem}\label{SettingirreducibleOrbits} Let $M$ be a submonoid of $\grAdj$, and $x\in\Aa{M}$. Let $D_1$ be an irreducible subgroup
of $\Ae{M^*}$, let $D_2$ be an irreducible subgroup  of $\Ae{M}$. Then the $D_1\times D_2$-orbit of the element 
$\al(1,x)\in \Spm\FK{M}$ is irreducible.
\end{Theorem}
\Proof
Let $Or$ be the $D_1\times D_2$-orbit of $\al(1,x)$, i.e., $Or=\al(D_1,x D_2)=\al(x^* D_1, D_2)$. Let $A_1$ and $A_2$ be closed subsets of 
$\Spm\FK{M}$, such that $Or\subseteq A_1\cup A_2$. We have to show $Or\subseteq A_1$ or $Or\subseteq A_2$.\vspace*{1ex}\\
Let $d_1\in D_1$. Due to the last proposition the map $\gamma_{d_1}:D_2\to \Spm\FK{M}$ defined by $\gamma_{d_1}(d_2):=\al(x^* d_1, d_2)$, 
$d_2\in D_2$, is a morphism.\\ 
Similarly, for $d_2\in D_2$, the map $\delta_{d_2}:D_1\to \Spm\FK{M}$ defined by $\delta_{d_2}(d_1) := \al(d_1,x d_2)$, $d_1\in D_1$, is a 
morphism.\vspace*{1ex}\\
Let $d_1\in D_1$. Because of $\gamma_{d_1}(D_2)\subseteq Or$, we have $\gamma_{d_1}^{-1}(A_1)\cup \gamma_{d_1}^{-1}(A_2)=D_2$. Furthermore 
$\gamma_{d_1}^{-1}(A_1)$ and $\gamma_{d_1}^{-1}(A_2)$ are closed. Because of the 
irreducibility of $D_2$ we get $\gamma_{d_1}^{-1}(A_1)=D_2$ or $\gamma_{d_1}^{-1}(A_2)=D_2$.\vspace*{1ex}\\  
Therefore the sets
\begin{eqnarray*}
  B_1 &:=& \Mklz{d_1\in D_1}{\gamma_{d_1}^{-1}(A_1)=D_2}\;\;,\\
  B_2 &:=& \Mklz{d_1\in D_1}{\gamma_{d_1}^{-1}(A_2)=D_2}
\end{eqnarray*}
satisfy $B_1\cup B_2 =D_1$.\\
Note that for $d_1\in D_1$ and $d_2\in D_2$ we have $\gamma_{d_1}(d_2)=\delta_{d_2}(d_1)$. The set $B_1$ is closed, because of 
\begin{eqnarray*}
  B_1 &=& \Mklz{d_1\in D_1}{ \gamma_{d_1}(d_2)\in A_1 \;\mb{ for all }\; d_2\in D_2}\;\;=\;\;\bigcap_{d_2\in D_2}\underbrace{\delta_{d_2}^{-1}(A_1)}_{closed}\;\;.
\end{eqnarray*}
Similarly, the set $B_2$ is closed.
Because of the irreducibility of $D_1$ we get $B_1=D_1$ or $B_2=D_1$, which is equivalent to $Or\subseteq A_1$ or $Or\subseteq A_2$.\\
\End
We have $\Ae{\grAdj}=\grEnd$. Due to Proposition \ref{Setting2} we get a map 
$\al:\grEnd\times \grEnd\to\Spm\FK{\grAdj}$. The following technical proposition,
which follows immediately from the definition of the Zariski closures in $\grEnd$, is very useful for determining such closures:
\begin{Prop}\label{Setting3} Let $M$ be a submonoid of $\grAdj$. We have:
\begin{eqnarray*}
 \Ae{M}   &=&  \Mklz{\phi\in\grEnd}{\al(1,\phi) \mb{ factors to a homomorphism } \FK{M}\to\F }\;\;.\\
 \Ae{M^*} &=&  \Mklz{\phi\in\grEnd}{\al(\phi, 1) \mb{ factors to a homomorphism } \FK{M}\to\F }\;\;.\\
\end{eqnarray*}
\end{Prop}
The formal Kac-Moody group $G_f$ acts faithfully on $\bigoplus_{\La\in P^+}L(\La)$. We identify $G_f$ with the corresponding subgroup of $\grEnd$. 
Under this identification the minimal Kac-Moody group $G\subseteq G_f$ is identified with a subgroup of $\grAdj$, invariant under taking the adjoint.\\
The results of this section can be applied to the 
$\F$-valued points of the algebra of strongly regular functions $\FK{G}$, because the restriction of the 
coordinate ring $\FK{\grAdj}$ to $G$ coincides with $\FK{G}$.\\
Also the monoid $\GD$ can be identified with a submonoid of $\grAdj$. In \cite{M}, Theorem 5.14, we showed $\Aa{G}\,=\GD$. In this paper we determine 
$\Gq$. We show the bijectivity of the map $\ti{\al}\,:\,\Ae{G}\times_{\Aa{G}}\Ae{G} \to \Spm\FK{G}$.\\
To prove its surjectivity, we use an induction over $|J|$, $J\subseteq I$, describing $\Spm\FK{G_J}$. To prepare this proof,
in the next two sections we first determine the $\F$-valued points of some other coordinate rings.
%
%
%
%
%
%
%
%
\section{The $\F$-valued points of $\FK{T_J}$, $\FK{T}$, $\FK{T_{rest}}$ ($J\subseteq I$)}
%
%
%
%
Recall that the torus  $T$ of the Kac-Moody group can be described by the following isomorphism of 
groups:\vspace*{1.5ex}\\
 \hspace*{8em} $\begin{array}{ccc}
 H\otimes_\Z\F^\times &\to & \;\;\:T \\
 \sum_{i=1}^{2n-l}h_i\otimes s_i\;\, &\mapsto &\, \prod_{i=1}^{2n-l} t_{h_i}(s_i) 
\end{array}$\vspace*{1.5ex}\\
The group algebra $\FK{P}$ of the lattice $P$ can be identified with the classical coordinate ring on $T$,
identifying $\sum c_\la e_\la\in \FK{P}\,$ with the function on $T$, which is defined by
\begin{eqnarray*}
  \left(\sum_\la c_\la e_\la\right)\left(\prod_{i=1}^{2n-l}t_{h_i}(s_i)\right) \;\,:=\,\;\sum_\la \,c_\la \,
  \prod_{i=1}^{2n-l}(s_i)^{\la(h_i)} &\;\,,\;\,&   (s_i\in\F^\times )\;\;.
\end{eqnarray*}
We have similar descriptions for the tori $T_J:=\Mklz{\prod_{j\in J}t_{h_j}(s_j)}{s_j\in\F^\times }$, $T_{rest}$, and its classical coordinate rings, replacing the lattices $H$, $P$  
by $H_J$, $P_J:=\Z\mb{-span}\Mklz{\La_j}{j\in J} $ or $H_{rest}$, $P_{rest}:=\Z\mb{-span}\Mklz{\La_i}{i=n+1,\ldots,2n-l} $, 
($J\subseteq I$).\\\\ 
In \cite{M}, Proposition 5.1, we determined the coordinate rings $\FK{T_J}$, $\FK{T}$ and $\FK{T_{rest}}$. They are in 
general only subalgebras of the classical coordinate rings of these tori:
\begin{Theorem}\label{SpT1}
Let $J\subseteq I$. Let $p_J:\,P \to P_J$ be the projection defined by $p_J(\la):= \sum_{j\in J} \la(h_j)\La_j$. 
We have: \vspace*{0.5ex}\\
1) $\qquad \FK{T_J}\:\,=\,\:\FK{\,p_J\left(X\cap P\right)}\;$.\vspace*{0.5ex}\\
2) $\qquad \FK{T}\,\:=\,\:\FK{X\cap P}\;$.\vspace*{0.5ex}\\
3) $\qquad \FK{T_{rest}}\,\:=\,\:\FK{P_{rest}}\;$.
\end{Theorem}
In \cite{M}, Theorem 5.2, we described the relative closures of $T_J$, $T$ and $T_{rest}$ in $\grAdj$. In the same way, but replacing 
$\grAdj$ by $\grEnd$, we can determine the closures of $T_J$, $T$ and $T_{rest}$. The proof of \cite{M}, Theorem 5.2, also shows, that 
$\{1\}\times \TD_J$, $\{1\}\times \TD$, resp. $\{1\}\times T_{rest}$ map bijectively to $\Spm\FK{T_J}$, $\Spm\FK{T}$, resp. $\Spm\FK{T_{rest}}$. 
Therefore we get: 
\begin{Theorem}\label{SpT2}
Let $J\subseteq I$.\vspace*{0.5ex}\\
1) We have $\Ae{T_J}\:=\:\TD_J$, and $\{1\}\times \Ae{T_J}$ maps bijectively to $\Spm\FK{T_J}$.\vspace*{0.5ex}\\
2) We have $\Tq\:=\:\TD$, and $\{1\}\times\Tq$ maps bijectively to $\Spm\FK{T}$.\vspace*{0.5ex}\\
3) $T_{rest}\,$ is closed, and $\{1\}\times T_{rest}$ maps bijectively to $\Spm\FK{T_{rest}}$.
\end{Theorem}
%
%
%
%
\section{The $\F$-valued points of $\FK{U_J}$ and $\FK{U^J}$ ($J\subseteq I$)\label{SpU}}
%
%
%
The coordinate ring of $U$ has been described by Kac and Peterson in \cite{KP2}, Lemma 4.3. From this Lemma follows immediately
part 1) of the next theorem. Part 2) has been shown in \cite{M}, Theorem 5.6.
\begin{Theorem}\label{SpU1} Let $J\subseteq I$.\\ 
1) Let $\La\in F_{I\setminus J}\cap P^+$, and $v_\La\in L(\La)_\La\setminus\{0\}\,$. Then:\vspace*{0.5ex}\\
\hspace*{1em} $\FK{U_J}$ is a symmetric algebra in $\{\,f_{v_\La yv_\La}\res{U_J}\,|\,y\in \n_J^-\,\}\,$. 
\vspace*{0.5ex}\\
2) Let $N\in F_J\cap P^+$, and $v_N\in L(N)_N\setminus\{0\}\,$. Then:\vspace*{0.5ex}\\
\hspace*{1em} $\FK{U^J}$ is as algebra generated by $\{\,f_{v_N yv_N}\res{U^J}\,|\,y\in (\n^J)^-\,\}\,$.
\end{Theorem}
Using these descriptions of the coordinate rings, it is possible to determine its $\F$-valued points:
\begin{Theorem}\label{SpU2}
Let $J\subseteq I\,$.\vspace*{0.5ex}\\
1) We have $\Ae{U_J} = (U_f)_J$, and $\{1\}\times \Ae{U_J}$ maps bijectively to $\Spm\FK{U_J}$.\vspace*{0.5ex}\\
2) We have $\Ae{U^J} = (U_f)^J$, and $\{1\}\times \Ae{U^J}$ maps bijectively to $\Spm\FK{U^J}$.
\end{Theorem}
{\bf Proof of 1)}: The case $J=\emptyset$ is trivial. Let $J$ be nonempty.
We first show $(U_f)_J\subseteq\Ae{U_J}$. The coordinate ring $\FK{\grEnd}$ is a symmetric algebra in the linear span of 
the functions $f_{vw}\res{gr\mb{-}End}$, $v,w\in L(\La)$, $\La\in P^+$, which is isomorphic to $\bigoplus_{\La\in P^+}L(\La)\otimes L(\La)$. 
We get a representation of the Lie algebra $(\n_f)_J$, by assigning 
$x\in (\n_f)_J$ the derivation $\delta_x: \FK{\grEnd}\to\FK{\grEnd}$, which is defined by
\begin{eqnarray*}
  \delta_x(f_{vw}\res{gr\mb{-}End}) \;\, := \,\; f_{v\,xw}\res{gr\mb{-}End} &\quad,\quad & v,w\in L(\La)\;\,,\,\;\La\in P^+\;\;. 
\end{eqnarray*}
The derivation $\delta_x$ is locally nilpotent, because $x\in(\n_f)_J$ acts locally nilpotent on $L(\La)$, $\La\in P^+$. The
homomorphism of algebras $exp(\delta_x)$ satisfies
\begin{eqnarray*}
   \al(1,1)\;\circ\; exp(\delta_x) &=& \al(1,exp(x))\;\;.
\end{eqnarray*}
Let $I(U_J)$ be the vanishing ideal of $U_J$ in $\FK{\grEnd}$. Due to the last identity and Proposition \ref{Setting3},
it is sufficient to show $\delta_x \left( I(U_J)\right) \subseteq I(U_J)$ for all $x\in (\n_f)_J$. \\
Let $f\in I(U_J)$. Let $x\in\g_\al$, $\al\in (\W_J)^+_{re}$. For all $u\in U_J$ and $t\in\F$ we have
\begin{eqnarray*}
   0\;=\; f(u\exp(tx)) \;=\;f(u)\,+\,t \delta_x(f)(u)\,+\,O(t^2)\;\;.
\end{eqnarray*}
The right side is polynomial in $t$. Due to $|\,\F\,|=\infty$ the coefficients of the powers of $t$ vanish, and we find 
$\delta_x(f)\in I(U_J)$.\\
The elements $x\in\g_\al$, $\al\in (\W_J)^+_{re}$, generate $\n_J$. Therefore also $\delta_x(f)\in I(U_J)$ for all $x\in\n_J$.\\
For an element $v\in L(\La)$, $\La\in P^+$, there exist only finitely many roots $\al\in\pW$ with $\g_\al v\neq\{0\}$. Therefore
for an element $x\in(\n_f)_J$ there exists an element $\ti{x}\in\n_J$, depending on $f$, such that 
$\delta_x(f)=\delta_{\ti{x}}(f)\in I(U_J)$. \vspace*{1ex}\\
We have $(U_f)_J\subseteq \Ae{U_J}$, and $\{1\}\times \Ae{U_J}$ maps injectively to $\Spm\FK{U_J}$. To prove 1), it remains
to show that $\{1\}\times (U_f)_J$ maps surjectively to $\Spm\FK{U_J}$.\\
Let $\La\in F_{I\setminus J}\cap P^+$, and $v_\La\in L(\La)_{\La}\setminus\{0\}$. Due to the description of $\FK{U_J}$ of the last theorem, 
it is sufficient to show, that for every linear map 
\begin{eqnarray*}
   l:\;\Mklz{f_{v_\La\, yv_{\La}}\res{U_J}}{y\in \n_J^-} &\to &\F
\end{eqnarray*}
there exists an element $\phi\in (U_f)_J$ with $l( f_{v_\La\,yv_\La}\res{U_J})=\al(1,\phi)(f_{v_\La\,y v_\La}\res{U_J})$ 
for all $y\in\n_J^-$.\\
Choose $\iBl$-dual bases of $\n_J^+$, $\n_J^-$, adopted to the root space decomposition:
\begin{eqnarray*}
\begin{array}{c}
  x_{\al i}\in\g_\al\quad\al\in\W_J^+\;,\;i=1,\ldots,m_\al \\
  y_{\beta i}\in\g_{-\beta}\quad\beta\in\W_J^+\;,\;i=1,\ldots,m_\beta \end{array} 
&\mb{ such that } & \iB{x_{\al i}}{y_{\beta j}}\,=\,\delta_{\alpha \beta}\delta_{ij}\;\;.
\end{eqnarray*}
We have $\iB{\La}{\beta}>0$ for all $\beta\in\W_J^+$. Define recursively elements $b_{\beta j}\in\F$, $j=1,\ldots,m_\beta$, $\beta\in\W^+_J$, 
as follows:
\begin{eqnarray*}
   b_{\beta 1} \;\,:=\,\; \frac{1}{(\,\La\,|\,\beta\,)\kB{v_\La}{v_\La}}\, l(f_{v_\La\,y_{\beta 1} v_\La}\res{U_J})&\mb{ for }& \Ht\,\beta=1\;\;.
\end{eqnarray*}
Let $\beta\in\W_J^+$ with $\Ht\,\beta>1$, and let $b_{\al i}$ be defined for all $i=1,\ldots, m_\al$, $\al\in\W_J^+$ with $\Ht\al<\Ht\beta$. Set:
\begin{eqnarray*}
   b_{\beta j} \,\;:=\;\, \frac{1}{\iB{\La}{\beta}\kB{v_\La}{v_\La}} \left( l(f_{v_\La\,y_{\beta j} v_\La}\res{U_J}) \,-\,
    \kkB{v_\La}{exp(\sum_{\al,i\atop ht\,\al<ht\,\beta}b_{\al i}x_{\al i}\,)\,y_{\beta j}v_\La} \right)\;\;.
\end{eqnarray*}
Next we show that we can take $\phi=\exp(\sum_{\al i}b_{\al i}x_{\al i})$. Let $\beta\in\W_J^+$. By expanding the exponential 
function we find for $j=1,\ldots,m_\beta$:
\begin{eqnarray*}
\lefteqn{ \kkB{v_\La}{exp(\sum_{\al,i}b_{\al i}x_{\al i}\,)\,y_{\beta j}v_\La}\;\,=}\\
 && \kkB{v_\La}{exp(\sum_{\al,i\atop ht\,\al\,<\,ht\,\beta}b_{\al i}x_{\al i}\,)\,y_{\beta j}v_\La} + 
    \kkB{v_\La}{exp(\sum_{\al,i\atop ht\,\al\,=\,ht\,\beta}b_{\al i}x_{\al i}\,)\,y_{\beta j}v_\La}\;\;.
\end{eqnarray*}
Denote by $\nu:\h\to\h^*$ the linear isomorphism induced by the invariant bilinear form $\iBl$.
Using $[x_{\beta i}, y_{\beta j}]=\delta_{ij}\,\nu^{-1}(\beta)$, we find that the second summand on the right side equals 
$b_{\beta j} \iB{\La}{\beta}\kB{v_\La}{v_\La}$. After inserting the definition of $b_{\beta j}$, the right side equals 
$l(f_{v_\La\,y_{\beta j} v_\La}\res{U_J})$.\\
\End
To prepare the proof of 2), we first show two propositions:
\begin{Prop} \label{SpU3} Let $\La\in P^+$ and $J\subseteq I$. Then $(U_f)^J$ fixes the points of $U(\n_J^-)L(\La)_\La$.
\end{Prop}
\Proof
The case $J=\emptyset$ is obvious, let $J\neq\emptyset$. If we fix an element $v\in U(\n_J^-)L(\La)_\La$, then for every 
element $\ti{x}\in (\n_f)^J$, there exists an element $x\in \n^J$, such that 
$\exp(\ti{x})v=\exp(x)v$. Therefore it is sufficient to show, that every element $x\in \n^J$ acts trivially on 
$U(\n_J^-)L(\La)_\La$.\\
Let $v_\La\in L(\La)_\La\setminus\{0\}$. We show by induction over $n\in\Nn$:\vspace*{1ex}\\
\hspace*{0.2em} $n=0\,: \quad\qquad\qquad x\, v_\La \:=\: 0 \qquad \mb{for all }\,x\in \n^J\;$.\vspace*{1ex}\\ 
\hspace*{0.2em} $n\in\N\,:\quad\; x\, y_1\cdots y_n v_\La\:=\: 0 \qquad \mb{for all }\,x\in \n^J,\;y_1,\,\ldots,\, y_n\in \n_J^-\;$.
\vspace*{1ex}\\
Clearly the statement for $n=0$ is valid. The induction step from $n$ to $n+1$ follows from the equation
\begin{eqnarray*}
   x\, y_1\,y_2\cdots y_{n+1} v_\La &=& [x\,,\, y_1]\,y_2\cdots y_{n+1} v_\La \;+\; y_1\,\left( x\,y_2\cdots y_{n+1} v_\La \right)\;,
\end{eqnarray*}
together with $[\n^J,\n_J^-]\subseteq \n^J$.\\
\End
\begin{Prop} \label{SpU4} Let $J\subseteq I$, and let $\La\in F_{I\setminus J}\cap P^+$, $N\in F_J\cap P^+$. The comorphism 
$m^*:\FK{U}\to \FK{U_J}\otimes\FK{U^J}$, dual to the multiplication map $m:U_J\times U^J\to U$, is an isomorphism of algebras.
Furthermore we have:
\begin{eqnarray}
   m^*(f_{v_\La\,yv_\La}\res{U})\;=\; f_{v_\La\,yv_\La}\res{U_J}\,\otimes\: 1   &\quad,\quad &  y\in\n_J^- \label{U-Formel1}\;\;.\\
   m^*(f_{v_N\,yv_N}\res{U})\;=\; 1\,\otimes\, f_{v_N\,yv_N}\res{U^J}   &\quad,\quad &  y\in(\n^J)^- \label{U-Formel2}\;\;.
\end{eqnarray}
\end{Prop}
\Proof 
Due to \cite{KP2}, Lemma 4.2, $\FK{U}$ is a Hopf algebra. Therefore for $f\in\FK{U}$ there exist functions $f_i,g_i\in\FK{U}$, 
$i=1,\ldots,m$, such that $f(u_1u_2)=\sum_{i=1}^m f_i(u_1)g_i(u_2)$ for all $u_1,u_2\in U$. Restricting to 
$u_1\in U_J$, $u_2\in U^J$ we get $f\circ m =\sum_{i=1}^{m} f_i\res{U_J}\otimes\: g_i\res{U^J}\,\in\,\FK{U_J}\otimes\FK{U^J}$.\\
$m^*$ is injective, because $m$ is surjective. To show the surjectivity of $m^*$, due to Theorem \ref{SpU1}, it is sufficient to 
show the equations (\ref{U-Formel1}) and (\ref{U-Formel2}). Due to the last proposition  we find for all $y\in \n_J^-$, and for 
all $u_1\in U_J$, $u_2\in U^J\subseteq (U_f)^J$:
\begin{eqnarray*}
  m^*(f_{v_\La\,yv_\La}\res{U})(u_1,u_2) &=& \kB{v_\La}{u_1u_2 y v_\La}\;\,=\;\,  \kB{v_\La}{u_1 y v_\La} \\
                                         &=& \left(f_{v_\La\,yv_\La}\res{U_J}\otimes\,1\,\right)(u_1,u_2)\;\;.
\end{eqnarray*}
For $j\in J$ the $\al_j$-string of $P(N)$ through the highest weight $N$ consists only of $N$, due to $N=\sigma_j N$. Therefore $G_J$ is 
contained in the stabilizer of $v_N$, and we find for $y\in (\n^J)^-$, and for all $u_1\in U_J$, $u_2\in U^J$:
\begin{eqnarray*}
  m^*(f_{v_N\,yv_N}\res{U})(u_1,u_2) &=& \kB{(u_1)^*v_N}{u_2 y v_N}\;\,=\;\, \kB{v_N}{u_2 y v_N}  \\
                                     &=& \left(\,1\,\otimes\, f_{v_N\,yv_N}\res{U^J}\right)(u_1,u_2)\;\;.
\end{eqnarray*}
\Ende
{\bf Proof of Theorem \ref{SpU2}, 2):} The case $J=I$ is trivial. Let $J\neq I$. We first show $(U_f)^J\subseteq \Ae{U^J}$:
Let $\phi\in (U_f)^J$. We have $(U_f)^J\subseteq U_f= \Ae{U}\subseteq \Ae{G}$. Denote by $\al(1,\phi):\FK{G}\to \F$ the homomorphism of algebras 
corresponding to $\phi$. Denote by $I(U^J)$ the vanishing ideal of $U^J$ in $\FK{G}$. Because of Proposition \ref{Setting3} it is sufficient 
to show that $I(U^J)$ is contained in the kernel of $\al(1,\phi)$. \\
Let $f\in I(U^J)$. Due to the Peter and Weyl theorem $f$ is of the form $f=\sum_{i} f_{v_i w_i}\res{G}$. Because $\phi$ is of the form 
$\phi=exp(\prod_\al x_\al)$, $x_\al\in (\n^J)_\al$, $\al\in (\W^J)^+$, it is sufficient to show:
\begin{eqnarray}
           0   &=&   \sum_i \kB{v_i}{w_i}\;\;,\label{Nullsumme1}\\
           0   &=&   \sum_i \kB{v_i}{x_1 \cdots x_k w_i} \qquad \mb{ for all } x_1,\,\ldots,\, x_k\in \n^J,\;k\in \N\;\;.\label{Nullsumme2}   
\end{eqnarray}
Equation (\ref{Nullsumme1}) follows because of $1\in U^J$. It is sufficient to show equation (\ref{Nullsumme2}) for a system of generators of $\n^J$.\\
Define recursively a multibracket for elements of $\g$:
\begin{eqnarray*}
  [x] &:=& x  \quad,\quad x \in \g \;\;,\\ 
 \mb{} [ x_{k+1},\, x_k,\,\ldots,\, x_1 ]   &:=& [x_{k+1},[x_k,\,\ldots,\, x_1]]  \quad,\quad x_1,\ldots ,x_k\in \g,\;k\in\N \;\;.
\end{eqnarray*}
By an easy induction, a multibracket of the form $[y,\,x_1,\,\ldots\,x_k]$ is a linear combination of multibrackets $[x_{i_1},\ldots,x_{i_k},y]$, where 
$(i_1,\ldots,i_k)$ is a permutation of $(1,\ldots,k)$.\\ 
$\n^+$ is generated by $\g_\al$, $\al\in\prW=(\W_J)^+_{re}\,\dot{\cup}\,(\W^J)^+_{re}$. We have $\n=\n_J\oplus \n^J$ and $[\n_J,\n^J]\subseteq \n^J$. All 
multibrackets with all elements in $\g_\al$, $\al\in(\W_J)^+_{re}$ are in $\n_J$, all multibrackets with an element in $\g_\beta$, $\beta\in (\W^J)^+_{re}$ 
are in $\n^J$. Therefore $\n^J$ is generated by all multibrackets, which contain an element in $\g_\beta$, $\beta\in (\W^J)^+_{re}$. By the remark of above 
it is sufficient to consider the multibrackets 
\begin{eqnarray*}
      [x_{\gamma_m},\,\cdots,\, x_{\gamma_1},\, x_{\gamma_0}] 
\end{eqnarray*}
with $x_{\gamma_i}\in \g_{\gamma_i}$, $i=0,\ldots,m$, and $\gamma_1,\,\ldots,\,\gamma_m\in \prW$, $\gamma_0\in (\W^J)^+_{re}$, $m\in\N$.\\
$U^J$ is normal in $U$. Therefore for $t_0,\,\ldots,\, t_m\in \F$ we have
\begin{eqnarray*}
   \lefteqn{u(t_m,\ldots,t_1,t_0)\;\,:=\;\,}\\
  && \exp(t_m x_{\gamma_m})\cdots \exp(t_1 x_{\gamma_1})exp(t_0 x_{\gamma_0})\exp(-t_1 x_{\gamma_1})\cdots \exp(-t_m x_{\gamma_m})
   \;\,\in\;\,U^J\;\;. 
\end{eqnarray*}
We consider a product $p$ of $k$ factors of such expressions in different variables. To simplify our notation, we only write down a factor in the middle 
of $p$: 
\begin{eqnarray*}
  p &=& \cdots\,u(t_m,\ldots,t_1,t_0) \,\cdots  \quad \mb{ where }\quad \ldots,\, t_m,\,\ldots,\,t_1,\,t_0,\,\ldots \in \F\;\;.
\end{eqnarray*}
Because of $p\in U^J$ we get
\begin{eqnarray*}
   0 \;\,=\;\, f(p)\;\,=\;\,\sum_i \kB{v_i}{\,\cdots \,u(t_m,\,\ldots,\,t_1,\,t_0)\,\cdots \,w_i}\;\;. 
\end{eqnarray*}
Because the root vectors belonging to real roots act locally nilpotent, the right side is polynomial in $\cdots\,t_0,t_1,\ldots, t_m\,\cdots \,$. Because of 
$|\F|=\infty$ the coefficients of the monomials are zero. In particular for the monomial $\cdots\,t_m\cdots t_1 t_0\,\cdots$ we find
\begin{eqnarray*}
       \sum_i \kB{v_i}{\,\cdots\,[x_{\gamma_m},\,\ldots,\,x_{\gamma_1},\, x_{\gamma_0}]\,\cdots\, w_i}\;\,=\;\,0\;\;.
\end{eqnarray*}
Now we show $(U_f)^J\supseteq \Ae{U^J}$: Let $u\in \Ae{U^J}$. Denote by $\al$ the map 
$\Ae{U^-}\times\Ae{U} \to \Spm\FK{U}$, and by $\ti{\al}$ the maps 
$\Ae{U_J^-}\times\Ae{U_J} \to \Spm\FK{U_J}$, $\Ae{(U^J)^-}\times\Ae{U^J} \to \Spm\FK{U^J}$. Using 
the last proposition and its notation, and 1), there exists an element $\ti{u}\in U_f$ such that
\begin{eqnarray}\label{altiutialu}
   \al(1, \ti{u}) &=& \left(\ti{\al}(1,1)\,\otimes\,\ti{\al}(1,u) \right)\,\circ\,m^*\;\;.
\end{eqnarray}
Write $\ti{u}$ in the form $u_J u^J$ with $u_J\in (U_f)_J$, $u^J\in (U_f)^J$. Choose elements $\La\in F_{I\setminus J}\cap P^+$ and 
$v_\La\in L(\La)_\La\setminus\{0\}$, and apply the left and right sides of the last equation to the elements $f_{v_\La\, yv_\La}\res{U}$, $y\in\n_J^-$. 
Using equation (\ref{U-Formel1}) of the last proposition we get
\begin{eqnarray*}
  \kB{v_\La}{u_J u^J y v_\La} \;\,=\;\, \kB{v_\La}{y v_\La} &\mb{ for all }& y\in\n_J^-\;\;.
\end{eqnarray*}
Due to Proposition \ref{SpU3} the left side is equal to $\kB{v_\La}{u_J y v_\La}$. Using Theorem \ref{SpU1}, 1), and 1) 
we conclude $u_J=1$.\vspace*{1ex}\\ 
Insert $\ti{u}=u^J$ in equation (\ref{altiutialu}). Choose elements $N\in F_J\cap P^+$ and $v_N\in L(N)_N\setminus\{0\}$, and apply the 
left and right sides of (\ref{altiutialu}) to the elements $f_{v_N\, yv_N}\res{U}$, $y\in(\n^J)^-$. Using equation (\ref{U-Formel2}) of the last 
proposition we find
\begin{eqnarray*}
  \kB{v_N}{u^J y v_N} \;\,=\;\, \kB{v_N}{ u y v_N} &\mb{ for all }& y\in(\n^J)^-\;\;.
\end{eqnarray*}
We have $1,\,u^J\in (U_f)^J\subseteq \Ae{U^J}$. Using Theorem \ref{SpU1}, 2) we conclude $\ti{\al}(1,u^J)=\ti{\al}(1,u)$, 
from which follows $u=u^J\in (U_f)^J$.\\ 
\End
%
%
%
\section{The $\F$-valued points of $\FK{G}$ and the Birkhoff decomposition}
%
%
%
%
Let $J\subseteq I$. Define the following subsets of $\grEnd$:
\begin{eqnarray*}\label{Bcovering}
   \widehat{G_f} &:=& \GD \,U_f\;\,=\;\, \bigcup_{\hat{n}\in\hat{N}} U^\pm \hat{n} U_f \;\;,\\
    (\widehat{G_f})_J &:=& \GD_J\, (U_f)_J\;\,=\;\, \bigcup_{\hat{n}\in\hat{N}_J} U_J^\pm \hat{n} (U_f)_J \;\;.
\end{eqnarray*}
We call the unions on the right the Bruhat and Birkhoff coverings. During our investigation of $\Spm\FK{G}$, 
we will find, that $\widehat{G_f}$ is a monoid, and the Bruhat and Birkhoff coverings are really decompositions. 
Similar things hold for $(\widehat{G_f})_J$.\vspace*{1ex}\\
The monoid $\Gq$ contains $G$, as well as the closures $\Tq=\TD$ and $\Ae{U}=U_f$. Therefore it also contains 
$\widehat{G_f}$, and due to Proposition \ref{Setting2} we have:
\begin{eqnarray}\label{ESpmG}
   \widehat{G_f}\times\widehat{G_f} &\mb{ maps to } & \Spm\FK{G}\;\;.
\end{eqnarray} 
Similarly we get:
\begin{eqnarray}\label{ESpmGJ}
   (\widehat{G_f})_J\times(\widehat{G_f})_J &\mb{ maps to } & \Spm\FK{G_J}\;\;.
\end{eqnarray}
Our first aim is to show the surjectivity of (\ref{ESpmG}). The key step of the proof is an induction over $|J|$, showing 
the surjectivity of (\ref{ESpmGJ}). The next two theorems prepare the induction step. They relate $\Spm\FK{G_J}$ 
to $\Spm\FK{G_{J\setminus\{j\}}}$, $j\in J$.\vspace*{1ex}\\
For $\La\in P^+$ choose a nonzero  element $v_\La\in L(\La)_\La$ and define 
\begin{eqnarray}
    \gt_\La &:=& \frac{f_{v_\La v_\La}\res{G}}{\kB{v_\La}{v_\La}} \;\,\in\;\, \FK{G}\;\;.
\end{eqnarray} 
The function $\gt_\La$ is independent of the chosen element $v_\La$, and the chosen nondegenerate contravariant symmetric bilinear form 
on $L(\La)$. Set $\gt_i:=\gt_{\La_i}$, $i=1,\,\ldots,\, n$.\\
Kac and Peterson showed in \cite{KP2} by checking on the dense principal open set $U^-TU^+$ of $G$:
\begin{eqnarray}\label{multgt}
   \gt_\La \gt_{\La'}\;\,=\,\; \gt_{\La+\La'} \quad \mb{ for all }\quad\La,\La'\in P^+\;\;.
\end{eqnarray}
The next theorem gives a covering of $\Spm\FK{G_J}$ by principal open sets, which are build with these functions. A variant of this covering for the 
full spectrum of $\FK{G}$ has been given by Kashiwara in \cite{Kas}, Proposition 6.3.1.\\
Denote the action of $G_J\times G_J$ on $\FK{G_J}$, which is induced by the action $\pi$ of $G\times G$ on $\FK{G}$, also by $\pi$. Write $\gt_\La$ 
instead of $\gt_\La\res{G_J}$ for short.
\begin{Theorem}\label{Ueb1}
Let $\emptyset\neq J\subseteq I$. We have
\begin{eqnarray*}
  \lefteqn{ \bigcup_{g,h\in G_J}\; \;\bigcup_{j\in J} \;\; D_{Specm\,\F\,[G_J]}\,(\,\pi(g,h)\gt_j\,)}\\
 &=& \left\{ \begin{array}{ccl}
   \Spm\,\FK{G_J} &\mb{ if }& J \mb{ is not special }.\\
   \left(\,\Spm\,\FK{G_J}\,\right)\setminus\{\,\al(1,e(R(J)))\,\} &\mb{ if }& J \mb{ is special }.
     \end{array} \right.
\end{eqnarray*}
\end{Theorem}
\Proof
Suppose there exists an element $\al\in\Spm\FK{G_J}$ not contained in the union on the left. Then for all $g,h\in G_J$ and $j\in J$ we have 
$\al(\pi(g,h)\gt_j)=0$. Here $\al\circ\pi(g,h)$ is a homomorphism of algebras. Because of the multiplicative property (\ref{multgt}) 
we find $\al(\,\pi(g,h)\gt_\La\,)=0$ for all $g,h\in G_J$ and $\La\in P^+_J\setminus\{0\}$, where $P_J^+:=P_J\cap P^+= \Nn\mb{-span}\,\Mklz{\La_j}{j\in J}$.\\
Now $L_J(\La):=U(\n_J^-)L(\La)_\La$ is an irreducible $G_J$-module, $\La\in P_J^+$. Since $G_J L(\La)_\La$ spans 
$L_J(\La)$, we find $\al(f_{vw}\res{G_J})=0$ for all $v,w\in L_J(\La)$, $\La\in P_J^+\setminus\{0\}$. Due to 
Theorem 5.12 of \cite{M} this is impossible if $J$ is not special, and $\al=\al(1,e(R(J)))$ if $J$ is special.\\
\End
The principal open subsets, which have been used in the covering of $\Spm\FK{G_J}$, can be obtained by 
\begin{eqnarray*}
    D_{Specm\,\F\,[G_J]}\,(\,\pi(g,h)\gt_j\,) &=& \left(\,\Spm\,\FK{\, D_{G_J}(\pi(g,h)\gt_j)\,}\,\right) \res{\F\,[G_J]}\;\;.
\end{eqnarray*} 
The next theorem gives a product decomposition of the principal open sets $D_{G_J}(\gt_j)$, $j\in J$. (Set $L=J\setminus\{j\}$ and 
$\La=\La_j$, $j\in J$). It can be proved in the same way as Theorem 5.11 a), and b) in \cite{M}. A decomposition 
of the coordinate rings analogous to the second part of b) has been given in \cite{Kas}, Lemma 5.3.4 and 5.3.5.
\begin{Theorem}\label{Ueb2}
Let $L\subsetneqq J\subseteq I$. Set $(U^\pm_J)^L:=\bigcap_{\sigma\,\in\,{\cal W}_L}\sigma U^\pm_J\sigma^{-1}$. For $\La\in P_J\cap F_L$ we 
have:\vspace*{0.5ex}\\
a) $\,D_{G_J}(\gt_\La)\:=\:(U^-_J)^L\,G_L\, T_{J\setminus L}\,(U^+_J)^L\,$.\vspace*{1ex}\\
b) The multiplication map
\begin{eqnarray*}
 m:\;\; (U^-_J)^L\times G_L\times T_{J\setminus L}\times (U^+_J)^L &\to & D_{G_J}(\gt_\La) 
\end{eqnarray*}
\hspace*{0.8em} is bijective, and its comorphism 
\begin{eqnarray*}
 m^*:\;\; \FK{D_{G_J}(\gt_\La)} &\to & 
  \FK{(U^-_J)^L}\otimes\FK{G_L}\otimes \FK{P_{J\setminus L}}\otimes\FK{(U^+_J)^L}
\end{eqnarray*}
\hspace*{0.8em} exists, and is an isomorphism of algebras.
\end{Theorem}
\begin{Theorem}\label{SpG1} Let $J\subseteq I$. Then $(\widehat{G_f})_J\times (\widehat{G_f})_J$ maps surjectively to $\Spm\FK{G_J}$.
\end{Theorem}
\Proof
We show the surjectivity by induction over $|J|$. The case $J=\emptyset$ is trivial. For $J=\{j\}$ the map 
$\{1\}\times G_{\{j\}}\to\Spm\FK{G_{\{j\}}}$ is already surjective, because $(G_{\{j\}},\FK{G_{\{j\}}})$ can be identified with 
$(SL(2,\F),\FK{SL(2,\F)})$.\vspace*{1ex}\\
Now the step of the induction from $|J|\leq m$ to $|J|=m+1$, ($1<m<|I|$):\\
Let $j\in J$. Due to the induction assumption we have a surjective map
\begin{eqnarray*} 
  \ti{\al}:\,(\widehat{G_f})_{J\setminus\{j\}}\times(\widehat{G_f})_{J\setminus\{j\}} &\to& \Spm\FK{G_{J\setminus\{j\}}}\;\;.
\end{eqnarray*}
Because $\FK{P_{\{j\}}}$ is the classical coordinate ring of the torus, we get a bijective map 
\begin{eqnarray*}
 \ti{\al}:\,\{1\}\times T_{\{j\}}&\to &\Spm\FK{P_{\{j\}}}\;\;.
\end{eqnarray*} 
$U_J^\pm$, $\FK{U_J^\pm}$, and $(U_f)_J$ can be identified with $U(A_J)^\pm$, $\FK{U(A_J)^\pm}$, and $U(A_J)_f$. Due to Theorem \ref{SpU2}, 
Proposition \ref{Setting2}, and Remark 3) after Proposition \ref{Setting2}, we have bijective maps 
\begin{eqnarray*}
  \ti{\al} : \,\{1\}\times ((U_f)_J)^{J\setminus\{j\}} &\to & \Spm\FK{(U_J)^{J\setminus\{j\}}} \;\;,\\ 
  \ti{\al} : \, ((U_f)_J)^{J\setminus\{j\}}\times \{1\} &\to & \Spm\FK{(U^-_J)^{J\setminus\{j\}}}\;\;.
\end{eqnarray*}
Due to the last theorem an element $\beta\in\Spm\FK{D_{G_J}(\gt_j)}$ can be written in the form
\begin{eqnarray*}
  \beta &=& \left(\,\ti{\al}(u,1)\otimes \ti{\al}(x,y)\otimes \ti{\al}(1,t)\otimes \ti{\al}(1,\ti{u})\,\right)\circ m^*
\end{eqnarray*}
with $u,\ti{u}\in ((U_f)_J)^{J\setminus\{j\}}$, $t\in T_{\{j\}}$, and $x,y\in (\widehat{G_f})_{J\setminus\{j\}}$.\vspace*{1ex}\\
Let $N\in P^+$. Choose $\kBl$-dual bases of $L(N)$, by choosing $\kBl$-dual bases 
\begin{eqnarray*}
    (\,a_{\la i}\,)_{\,i=1,\ldots,m_\la} &\quad,\quad & (\,c_{\la i}\,)_{\,i=1,\ldots,m_\la} 
\end{eqnarray*}
of $L(N)_\la$ for every $\la\in P(N)$. Let $v\in L(N)_\la$, $w\in L(N)_\mu$. By applying $m^*$ to $f_{vw}\res{G_J}$ we find
\begin{eqnarray*}
 \lefteqn{ m^* (f_{vw}\res{G_J}) }\\ 
   &=&  \sum_{(\la',\,i)\,,\,\la'\geq \la \atop (\mu',\,j)\,,\,\mu'\geq \mu }
    f_{v a_{\la' i}}\res{(U_J^-)^{J\setminus\{j\}}} \otimes\, f_{c_{\la' i} a_{\mu' j}}\res{G_{J\setminus\{j\}}}\otimes \,
  e_{\mu'(h_j)\La_j}\res{T_{\{j\}}}\otimes\, f_{c_{\mu' j} w} \res{(U_J)^{J\setminus\{j\}}}\;\;.
\end{eqnarray*}
This sum is finite, due to $P(N)\subseteq N-Q_0^+$. By applying $\beta$ to $f_{vw}\res{G_J}$ we get
\begin{eqnarray*}
  \beta(f_{vw}\res{G_J}) &=& \sum_{(\la',\,i)\,,\,\la'\geq \la \atop (\mu',\,j)\,,\,\mu'\geq \mu } \kB{uv}{a_{\la' i}}
   \kB{x c_{\la' i}}{yt a_{\mu' j}} \kB{c_{\mu' j}}{\ti{u}w} \\
    &=& \kB{xuv}{yt\ti{u}w}\;\;.
\end{eqnarray*}
In particular for $v\in U(\n_J^-) L(N)_N \cap L(N)_\la$, $w\in U(\n_J^-) L(N)_N \cap L(N)_\mu$ we have
\begin{eqnarray*}
  \beta(f_{vw}\res{G_J}) \;\,=\;\,\kB{xuv}{yt\ti{u}w}\;\,=\;\, \al(xu,yt\ti{u})(f_{vw}\res{G_J})\;\;.
\end{eqnarray*}
Here $\al$ denotes the map $(\widehat{G_f})_J\times(\widehat{G_f})_J \to \Spm\FK{G_J}$ due to (\ref{ESpmGJ}).
Because $\FK{G_J}$ is already spanned by the functions $f_{vw}\res{G_J}$, $v,w\in U(\n_J^-)L(N)_N$, $N\in P_J^+=P_J\cap P^+$, we conclude 
$\beta\res{\F\,[G_J]}=\al(xu,yt\ti{u})\in\Spm\FK{G_J}$.\vspace*{1ex}\\
It is easy to check that for $g,h\in G_J$ we have
\begin{eqnarray*}
  \left(\,\Spm\,\FK{D_{G_J}(\pi(g,h)\gt_j)}\,\right)\res{\F\,[G_J]} &=&  
                                      \left(\,\Spm\,\FK{D_{G_J}(\gt_j)}\,\right)\res{\F\,[G_J]}\circ \:\pi(g^{-1},h^{-1})\;\;.
\end{eqnarray*}
Because of Theorem \ref{Ueb1}, and because of the equation $\al(a,b)\circ\pi(c,d)=\al(ac,bd)$, we get
\begin{eqnarray*}
  \lefteqn{ \bigcup_{g,h\in G_J} \,\bigcup_{j\in J} \,\Mklz{\al(xug^{-1},yt\ti{u}h^{-1}) }{ x,y\in (\widehat{G_f})_{J\setminus\{j\}}\,,
    \;u,\ti{u}\in ((U_f)_J)^{J\setminus\{j\}}\,,\;t\in T_{\{j\}} } }\\
  &=& \left\{ \begin{array}{ccl}
   \Spm\FK{G_J} &\mb{ if }& J \mb{ is not special }.\\
   \Spm\FK{G_J}\setminus\{\,\al(1,e(R(J)))\,\} &\mb{ if }& J \mb{ is special }.
     \end{array} \right.\qquad\qquad\qquad\qquad\qquad\qquad
\end{eqnarray*}
Using the Bruhat covering of $(\widehat{G_f})_J$, the Bruhat decompositions of $(G_f)_J$ and $\GD_J$, we find
\begin{eqnarray*}
 (\widehat{G_f})_J &\subseteq & (\widehat{G_f})_J(G_f)_J \,\;=\,\; U_J \ND_J (U_f)_J\, (G_f)_J \;\,=\;\,
  U_J\ND_J \, U_J N_J(U_f)_J\\
    &=& \GD_J (U_f)_J\;\,=\;\,(\widehat{G_f})_J\;\;.
\end{eqnarray*}
From this follows $(\widehat{G_f})_{J\setminus\{j\}} T_{\{j\}}((U_f)_J)^{J\setminus\{j\}} G_J \subseteq (\widehat{G_f})_J$ for all $j\in J$. 
Furthermore, if $J$ is special, then $e(R(J))\in  (\widehat{G_f})_J$. Therefore the map $\al: (\widehat{G_f})_J\times (\widehat{G_f})_J\to
\Spm\FK{G_J}$ is surjective.\\
\End
\begin{Cor}\label{SpG2} $\widehat{G_f}\times \widehat{G_f}$ maps surjectively to $\Spm\FK{G}$.
\end{Cor}
\Proof
Let $m: G_I\times T_{rest}\to G$ be the multiplication map, and $m^ *:\FK{G}\to\FK{G_I}\otimes \FK{T_{rest}}$ its comorphism. Due to the
bijectivity of $m^*$, and due the last theorem and Theorem \ref{SpT2}, 3), the elements of $\Spm\FK{G}$ are given by
\begin{eqnarray*}\label{tensorhom}
   \left(\ti{\al}(x,y)\otimes\ti{\al}(1,t)\right)\circ m^*\;\;,
\end{eqnarray*}
where $x,y\in (\widehat{G_f})_I$ and $t\in T_{rest}$. Applying this expression to the matrix coefficients $f_{vw}\res{G}$, 
$v\in L(\La)_\la$, $w\in L(\La)_\mu$, $\la,\mu\in P(\La)$, $\La\in P^+$, we find
\begin{eqnarray*}
  \left(\,\ti{\al}(x,y)\otimes\ti{\al}(1,t)\,\right) (m^*(f_{vw}\res{G})) &=& \left(\ti{\al}(x,y)\otimes\ti{\al}(1,t)\right)
   (f_{vw}\res{G_I}\otimes\, e_\mu\res{T_{rest}})\\ 
   &=&\kB{xv}{yw}e_\mu(t) \;\,=\,\; \al(x,yt)(f_{vw}\res{G})\;\;.
\end{eqnarray*} 
Here $\al$ denotes the map $ \widehat{G_f}\times\widehat{G_f} \to \Spm\FK{G}$ due to (\ref{ESpmG}).
Therefore $\Spm\FK{G}=\Mklz{\al(x,\ti{y})}{x\in(\widehat{G_f})_I\,,\,\ti{y}\in (\widehat{G_f})_I T_{rest} }$, in particular $\al$ is surjective.\\
\End
The next theorem gives one of the main results of this paper: A description of the $\widehat{G_f}\times\widehat{G_f}$-set $\Spm\FK{G}$.
\begin{Theorem}\label{SpG3}\mb{}\vspace*{0.5ex}\\
1) We have $\Aa{G}\,=\,\widehat{G}$ and $\Ae{G}\,=\,\widehat{G_f}$. In particular $\widehat{G_f}$ 
is a monoid.\vspace*{0.5ex}\\ 
2) The map $\Gq\times_{\Aa{G}} \Gq\to \Spm\FK{G}$ is a $\Gq\times\Gq$-equivariant bijection.
\end{Theorem}
\Proof
The first equation of 1) has been shown in  \cite{M}, Theorem 5.14. To prepare the proof of the rest of the theorem, we first show the
following statements a) - d):\vspace*{0.5ex}\\
a) We show $\Ae{B^-}\subseteq \,\Aa{B^-}$: 
For $\La\in P^+$ and $\la\in P(\La)$ define
\begin{eqnarray*}
    L(\La)_{\la\,\downarrow} &:=& \bigoplus_{\mu\,\in \,P(\La)\cap (\la-Q_0^+)} L(\La)_\mu\;\;.
\end{eqnarray*}
Since $L(\La)_{\la\,\downarrow}$ is $B^-$-invariant, we have
\begin{eqnarray*}
 \kB{w}{B^- v}\;=\; 0 &\mb{ for all}& v\in L(\La)_{\la\,\downarrow}\;\,,\;\,w\in \left(L(\La)_{\la\,\downarrow}\right)^\bot \;\;.
\end{eqnarray*}
These equations are also valid, if $B^-$ is replaced by its Zariski closure $\Ae{B^-}$, and due to the orthogonality 
of the weight spaces we have $\left(L(\La)_{\la\,\downarrow}\right)^{\bot\, \bot} = L(\La)_{\la\,\downarrow} $. Therefore 
$L(\La)_{\la\,\downarrow}$ is also $\Ae{B^-}$-invariant.\vspace*{0.5ex}\\ 
Let $b\in \Ae{B^-}$. Choose a pair of $\kBl$-dual bases of $L(\La)$, by choosing $\kBl$-dual bases
\begin{eqnarray*}
       (\,a_{\la i}^\La\,)_{\,i=1,\ldots,m_\la} &\quad,\quad &
       (\,c_{\la i}^\La\,)_{\,i=1,\ldots,m_\la} 
\end{eqnarray*}
of $L(\La)_\la$ for every $\la\in P(\La)$. For a fixed weight $\mu\in P(\La)$ we have 
\begin{eqnarray*}
  \kkB{a^\La_{\mu i}}{b\,c^\La_{\mu' i'}} & \neq&  0 
\end{eqnarray*}
at most for the finitely many weights $\mu'\in P(\La)$ with $\mu'\geq \mu$. Therefore we get a well defined linear map 
$\psi_\La$ by
\begin{eqnarray*}
   \psi_\La v \;\,:=\,\; \sum_{\mu\,i}\;\sum_{\mu'\,i'} \;a_{\mu'\, i'}^\La \kB{a^\La_{\mu i}}{b\,
  c^\La_{\mu' i'}} \kB{c^\La_{\mu i}}{v} &\;,\;& v\in L(\La)\;\;.
\end{eqnarray*}
These maps $\psi_\La$, $\La\in P^+$, define an element $\psi$ of $\grEnd$. It is easy to
check, that $\psi$ is the adjoint of $b$. Therefore $ b\in \Ae{B^-} \cap \grAdj = \Aa{B^-}$.\vspace*{1ex}\\
b) Let $\psi\in\Gq$, $b\in \Ae{B}$ with $\al(1,\psi)=\al(b,1)$. We show
\begin{eqnarray*}
  b,\,\psi\in\GD &\mb{and}& \psi^*\;=\;b \;\;.
\end{eqnarray*}
Due to the second part of Proposition \ref{Setting3}, the homomorphism 
$\al(b,1) =\al(1,\psi)$ factors to a homomorphism $\FK{B^-}\to\F$. Due to the first part of Proposition \ref{Setting3},  
and part a) of above, we find
\begin{eqnarray*} 
    \psi\;\,\in\;\,\Ae{B^-}\;\,\subseteq\;\,\Aa{B^-}\;\,\subseteq \;\,\Aa{G} 
    \,\;= \,\;\GD\,\;\subseteq \,\;\grAdj\;\;.
\end{eqnarray*}
By checking on the matrix coefficients, using the non-degeneracy of the contravariant bilinear forms, we find $\psi^*=b$.
\vspace*{0.5ex}\\
c) Let $g,h\in G_f$ and $x,y,\ti{x},\ti{y}\in \Gq$. By checking on the matrix coefficients we find:
\begin{eqnarray*}
  \al(xg,yh)\;=\;\al(\ti{x},\ti{y}) &\iff& \al(x,y)\;=\;\al(\ti{x}g^{-1},\ti{y}h^{-1})\;\;.
\end{eqnarray*}
d) We show $\,\Spm\FK{G} =\al(U_f,\ND U_f)\,$:
Due to the last corollary, every element of $\Spm\FK{G}$ is of the form $\al(\psi,\ti{\psi})$ with 
$\psi,\ti{\psi}\in\widehat{G_f}$. Due to the Birkhoff covering of $\widehat{G_f}$, we can write $\psi$, $\ti{\psi}$ in 
the form
\begin{eqnarray*}
\begin{array}{ccc}
  \psi & = & u_-n_\sigma e(R) u_+ \\
  \ti{\psi} & = & \ti{u}_- \ti{n}_{\ti{\sigma}} e(\ti{R}) \ti{u}_+
\end{array} &\;\mb{ with }\; & u_-,\,\ti{u}_-\in U^-\;\;,\;\; u_+,\,\ti{u}_+ \in U_f\;\;,\;\; n_\sigma,\,\ti{n}_{\ti{\sigma}}\in N\;\;.
\end{eqnarray*}
From this follows:
\begin{eqnarray*}
      (\psi,\ti{\psi})  & \sim & (u_+\, ,\, e(R)n_\sigma^* u_-^*\ti{u}_-\ti{n}_{\ti{\sigma}} e(\ti{R}) \ti{u}_+)\;\;.
\end{eqnarray*}
Due to the Birkhoff decomposition of $\GD$, we can write $e(R)n_\sigma^* u_-^*\ti{u}_-\ti{n}_{\ti{\sigma}} e(\ti{R})$
in the form $u'_- \hat{n}' u'_+$ with $\hat{n}'\in\ND$, $u'_\pm\in U^\pm$. We get
\begin{eqnarray*}
   (\psi,\ti{\psi})  & \sim & ((u_-')^*u_+\, ,\, \hat{n}' u_+'\ti{u}_+)\;\;.
\end{eqnarray*}
Therefore $\al(\psi,\ti{\psi})   = \al((u_-')^*u_+\, ,\, \hat{n}' u_+'\ti{u}_+)$.\vspace*{1ex}\\
Now we can prove the theorem:\vspace*{1ex}\\
To 1) To show $\Ae{G}\,=\,\widehat{G_f}$, it is sufficient to show $\Ae{G}\,\subseteq\,\widehat{G_f}$. Let 
$\phi\in \Ae{G}$. Then due to d) there exist elements $u,\ti{u}\in U_f$, $e(R)n_\sigma\in\ND$ such that
\begin{eqnarray*}
  \al(1,\phi)  & = & \al( u\, ,\, e(R)n_\sigma \ti{u})\;\;.
\end{eqnarray*}
Using c) this is equivalent to
\begin{eqnarray*}
   \al(1,\phi(n_\sigma\ti{u})^{-1}) &=& \al(e(R)u,1)\;\;.
\end{eqnarray*}
The monoid $\Ae{B}$ contains $\Tq$ and $U_f$. Therefore $e(R)u\in \Tq U_f\subseteq \Ae{B}$, and due to b) we get
$\phi\in \GD n_\sigma \ti{u}\subseteq \widehat{G_f}$.\vspace*{1ex}\\
To 2) Due to the last corollary, we only have to show the injectivity in 2). Due to the proof of d), we may start with elements
$u,\ti{u}\in U_f$, $g,\ti{g}\in N U_f$ such that 
\begin{eqnarray*}
  \al(e(R)u,g) &=&\al(e(\ti{R})\ti{u},\ti{g}) \;\,.
\end{eqnarray*}
Using c) this equation is equivalent to
\begin{eqnarray*}
  \al(e(R)u\ti{u}^{-1},1)  &=&  \al(1,e(\ti{R})\ti{g}g^{-1})\;\;.
\end{eqnarray*}
Due to  b) we find $e(R)u\ti{u}^{-1}\in\GD$ and $(e(R) u\ti{u}^{-1})^* \,=\, e(\ti{R})\ti{g}g^{-1}$. From this follows
\begin{eqnarray*}
   (\,e(R)u\,,\, g\,) &=& (\,(e(R)u\ti{u}^{-1})\ti{u}\,,\, g\,)\;\,\sim\;\, (\,\ti{u}\,,\, (e(R)u\ti{u}^{-1})^*g\,)
              \;\,=\;\,(\,\ti{u}\,,\, e(\ti{R})\ti{g}\,) \\
                     &\sim & (\,e(\ti{R})\ti{u}\,,\, \ti{g}\,) \;\;.
\end{eqnarray*}
\End
The group $U_f\times U_f$ acts on $\Spm\FK{G}$. The corresponding partition into orbits is described by the Birkhoff decomposition
in the following theorem:
\begin{Theorem}\mb{}\\
1) There are the following Bruhat- and Birkhoff decompositions of $\widehat{G_f}$:
\begin{eqnarray*}
   \widehat{G_f} &=& \dot{\bigcup_{\hat{n}\in\hat{N}}}\; U^\pm \hat{n} U_f\;\;.
\end{eqnarray*}
2) There is the following Birkhoff decomposition of $\Spm\FK{G}$:
\begin{eqnarray*}
   \Spm\FK{G} &=& \dot{\bigcup_{\hat{n}\in\hat{N}} }\; \al(U_f, \hat{n} U_f)\;\;.
\end{eqnarray*}
\end{Theorem}
\Proof
Due to the Bruhat and Birkhoff coverings of $\widehat{G_f}$ and part d) of the proof of the last theorem, we only have to 
show that these unions are disjoint.\\
a) First we do this for 2). Suppose there exist elements $u_1,u_2,\ti{u}_1,\ti{u}_2\in U_f$ such that
\begin{eqnarray*}
    \al(u_1, n_\sigma e(R) u_2)  &=&  \al(\ti{u}_1,\ti{n}_{\ti{\sigma}} e(\ti{R}) \ti{u}_2) \;\,.
\end{eqnarray*}
Using part c) of the proof of the last theorem, we find for all $v,w\in L(\La)$, $\La\in P^+$:
\begin{eqnarray}\label{kBGl1}
  \kB{v}{n_\sigma e(R) w} &=& \kkB{\ti{u}_1(u_1)^{-1}v}{\ti{n}_{\ti{\sigma}}e(\ti{R})\ti{u}_2 (u_2)^{-1}w}\;\; .
\end{eqnarray}
Fix an element $\La\in P^+$ with $P(\La)\cap X\setminus R\neq\emptyset$. Fix elements $\mu\in P(\La)\cap X\setminus R$ and 
$w_\mu\in L(\La)_\mu\setminus\{0\}$. By inserting $w=w_\mu$ in the last equation we find
\begin{eqnarray*}
   0 \;=\; \kkB{\ti{v}}{e(\ti{R})\ti{u}_2 (u_2)^{-1}w_\mu} &\mb{ for all }& \ti{v}\in L(\La)\;\;.
\end{eqnarray*}
Because of $\ti{u}_2(u_2)^{-1}\in U_f$ we find $\mu\in X\setminus \ti{R}$.\\
Since $\bigcup_{\La\in P^+}P(\La)=X\cap P$ this shows $(X\setminus R)\cap P\subseteq (X\setminus\ti{R})\cap P$, from which follows 
$R\supseteq \ti{R}$. We may interchange the variables with and without $\,\ti{\mb{}}\,$, and get $R=\ti{R}$.\vspace*{1ex}\\
Fix an element $\La\in P^+$ with $P(\La)\cap R\neq\emptyset$. Fix elements $\mu\in R\cap P(\La)$ and 
$w_\mu\in L(\La)_\mu\setminus\{0\}$. Because of $\ti{u}_2 (u_2)^{-1}\in U_f$ we have 
\begin{eqnarray*}
   \ti{\sigma}\mu   &\in &    supp(\ti{n}_{\ti{\sigma}} e(R) \ti{u}_2 (u_2)^{-1}w_\mu)\;\;.
\end{eqnarray*}
Choose a maximal weight $m$ of $supp(\ti{n}_{\ti{\sigma}} e(R) \ti{u}_2 (u_2)^{-1}w_\mu)$ with $m\geq \ti{\sigma}\mu$. 
By inserting $w=w_\mu$ in equation (\ref{kBGl1}), using $\ti{u}_1 (u_1)^{-1}\in U_f$, we find 
\begin{eqnarray}\label{kBGl3}
  \kB{v_m}{n_\sigma w_\mu} &=& \kB{v_m}{\ti{n}_{\ti{\sigma}}e(R)\ti{u}_2 (u_2)^{-1}w_\mu}\;\;\mb{ for all }\;v_m\in L(\La)_m\;\;.\;\;\;
\end{eqnarray}
Because $\kBl$ is nondegenerate on $L(\La)_m$, there exists an element $v_m$, such that the right side is nonzero. Therfore 
$\sigma \mu = m\geq \ti{\sigma}\mu$. Interchanging the variables with and without $\,\ti{\mb{}}\,$ gives $\ti{\sigma}\mu=\sigma\mu = m$.\\
Inserting $m=\ti{\sigma}\mu$ in (\ref{kBGl3}), and using $\ti{u}_2 (u_2)^{-1}\in U_f$, we find 
\begin{eqnarray*}
  \kB{v}{n_\sigma w_{\mu}} &=& \kB{v}{\ti{n}_{\ti{\sigma}}w_\mu} \;\;\mb{ for all }\;v\in L(\La)_m \,=\, L(\La)_{\sigma\mu} \,= 
  \,L(\La)_{\ti{\sigma}\mu} \;\;.
\end{eqnarray*}
Because of the orthogonality of the weight spaces, this equation is also valid for all $v\in L(\La)$.\\
This shows
\begin{eqnarray*}
  \kB{v}{n_\sigma e(R)w} \;=\; \kB{v}{\ti{n}_{\ti{\sigma}}e(R)w} &\mb{ for all }&
   v, w \in  L(\La)  \;,\;  \La\in P^+\;\;.
\end{eqnarray*}
Therefore we get $n_\sigma e(R) =\ti{n}_{\ti{\sigma}} e(R) = \ti{n}_{\ti{\sigma}} e(\ti{R})$. \vspace*{1ex}\\
b) The Birkhoff decomposition of $\widehat{G_f}$ follows from 2), by using the embedding of $\widehat{G_f}$ in $\Spm\FK{G}$. 
Suppose there exist elements $u_1,\ti{u}_1\in U^-$, $u_2,\ti{u}_2\in U_f$ such that
\begin{eqnarray*}
   (u_1)^*n_\sigma e(R) u_2  &=& (\ti{u}_1)^*\ti{n}_{\ti{\sigma}} e(\ti{R}) \ti{u}_2 \;\;.
\end{eqnarray*}
From this equation follows equation (\ref{kBGl1}), and as in part a) of the proof we can deduce $R=\ti{R}$. Also similar 
as in part a), but now using a minimal weight $m$ of $supp(\ti{n}_{\ti{\sigma}} e(R) \ti{u}_2 (u_2)^{-1}w_\mu)$ with 
$m\leq \ti{\sigma}\mu$, we find $n_\sigma e(R)= \ti{n}_{\ti{\sigma}}e(R)$.\\
\End
%
%
%
%
%
%
\section{The stratification of the spectrum of $\F$-valued points of $\FK{G}$ in $G_f\times G_f$-orbits}
%
%
%
%
%
In this section, we show that the $G_f\times G_f$-orbits of $\Spm\FK{G}$ are locally closed, irreducible, and in one to one correspondence 
with the finitely many special subsets of $I$. The closure relation is given by the inverse inclusion of the special sets. We give a 
countable covering of each orbit by big cells. We show that there exist stratified transversal slices to the orbits at any of their points.
\vspace*{1ex}\\
To cut short our notation, we denote by $x\tr y$ the image of $(x,y)\in\widehat{G_f}\times\widehat{G_f}$ under the surjective map 
$\al:\widehat{G_f}\times\widehat{G_f}\to\Spm\FK{G}$. We make use of this map as a parametrization of $\Spm\FK{G}$.\\
Recall that we have $x\tr zy = z^*x\tr y$, $x,y\in\widehat{G_f}$, $z\in \GD$. Recall that $\widehat{G_f}\times\widehat{G_f}$
acts on $\Spm\FK{G}$ by morphisms from the right, i.e.,
\begin{eqnarray*}
  (x\tr y) (\ti{x},\ti{y}) &=& x\ti{x}\tr y\ti{y} \quad,\quad  x,y,\ti{x},\ti{y}\in \widehat{G_f}\;.
\end{eqnarray*}
The Chevalley involution of $\FK{G}$ induces an involutive morphism $*$ on $\Spm\FK{G}$, which we also call Chevalley 
involution. It is given by the switch map:
\begin{eqnarray*}
  (x\tr y)^*  &=& y\tr x  \quad,\quad x,y\in \widehat{G_f}\;.
\end{eqnarray*}
In this section we also denote by $f_{vw}$ the function on $\Spm\FK{G}$, induced by the matrix coefficient $f_{vw}\res{G}\,\in\FK{G}$. It is given by
\begin{eqnarray*}
 f_{vw} (x\tr y)  &=& \kB{xv}{yw}  \quad,\quad x,y\in \widehat{G_f}\;.
\end{eqnarray*}
For a set $M\subseteq \Spm\FK{G}$ denote by $\Asp{M}$ its Zariski closure. Note, that for a subset $A\subseteq\GD$ we have $1\tr \Aa{A}\,\subseteq \Asp{1\tr A}$. Similarly, for
a subset $A\subseteq \widehat{G_f}$, we have $1\tr \Ae{A}\subseteq \Asp{1\tr A}$. These formulas are useful to determine closures in $\Spm\FK{G}$.
\begin{Theorem}\label{S1} 1) The partition of $\Spm \FK{G}$ in $G_f\times G_f$-orbits is given by
\begin{eqnarray}\label{partitionSpm}
   \Spm\FK{G} &=& \dot{\bigcup_{\Xi\;special}} \, G_f\tr e(R(\Xi))G_f\;\;.
\end{eqnarray}
2) Let $\Th$ be special. The orbit $G_f\tr e(R(\Th))G_f$ is locally closed and irreducible. Its closure is given by 
\begin{eqnarray}\label{closure}
   \dot{\bigcup_{\Xi\;special\;,\;\Xi\supseteq \Th}} G_f\tr e(R(\Xi))G_f\;\;.
\end{eqnarray}
\end{Theorem}
\Proof
a) We first decompose the $G_f\times G_f$-orbit $G_f\tr e(R(\Th))G_f$ in a union of $U_f\times U_f$-orbits, i.e., we show
\begin{eqnarray}\label{orbit}
   G_f\tr e(R(\Th))G_f &=& \dot{ \bigcup_{\hat{\sigma}\,\in\,{\cal W}\varepsilon(R(\Th)){\cal W}} } U_f\tr (\hat{\sigma} T) U_f\;\;.
\end{eqnarray}
By using the Birkhoff decomposition of $G_f$, and Proposition 2.14, Theorem 2.15 b) of \cite{M}, compare also the section preliminaries, we get
\begin{eqnarray*}
 G_f\tr e(R(\Th))G_f  &=& U^- N U_f \tr e(R(\Th))G_f \;\,=\;\,  U_f \tr N U e(R(\Th)) G_f  \\ &=&
   \bigcup_{\sigma\in\cal W} U_f \tr e(\sigma R(\Th)) G_f\;\;.
\end{eqnarray*}
Because $\We_{\Th\cup\Th^\bot}$ is the stabilizer of the face $R(\Th)$ as a whole, we may restrict the last union to the minimal 
coset representatives $\We^{\Th\cup\Th^\bot}$ of $\We/\We_{\Th\cup\Th^\bot}$, characterized by 
$\We^{\Th\cup\Th^\bot}=\Mklz{\sigma\in\We}{\sigma\al_i\in\prW \mb{ for all } i\in \Th\cup\Th^\bot}$. Next we
insert in $U_f \tr e(\sigma R(\Th)) G_f $ the Birkhoff decomposition $G_f=(\sigma U^-\sigma^{-1})NU_f$. By using Proposition 2.14 of 
\cite{M}, compare the section preliminaries, we get 
\begin{eqnarray*}
  \bigcup_{\sigma\in{\cal W}^{\Th\cup\Th^\bot}\,,\,\tau\in{\cal W} } U_f \tr 
   \underbrace{\sigma U^-_{\Th^\bot}\sigma^{-1}}_{\subseteq U^-} e(\sigma R(\Th)) \tau B_f  &=& 
  \bigcup_{\sigma\in{\cal W}^{\Th\cup\Th^\bot}\,,\,\tau\in{\cal W} } U_f \tr  e(\sigma R(\Th)) \tau B_f\;\;.
\end{eqnarray*}
Since $\bigcup_{\sigma\in{\cal W}^{\Th\cup\Th^\bot}\,,\,\tau\in{\cal W} } \ve{\sigma R(\Th)}\tau = \We \,\ve{R(\Th)}\,\We$, we have shown the 
equality in (\ref{orbit}). The disjointness of the union follows from the Birkhoff decomposition of $\Spm\FK{G}$.\\
Taking into account the the partition $\WeD = \dot{\bigcup}_{\Th \;sp}\, \We\,\ve{R(\Th)} \,\We$, part 1) 
of the theorem follows from (\ref{orbit}) the Birkhoff decomposition of $\Spm\FK{G}$.\vspace*{1ex}\\
b) Next we show that the union $\bigcup_{\Xi\;sp.\;,\;\Xi\supseteq \Th} G_f\tr e(R(\Xi))G_f$ is closed, i.e., we show that it 
is the common zero set of the functions
\begin{eqnarray*}
   f_{vw} \quad,\quad v,w\in L(\La)\quad,\quad \La\in P^+\setminus \overline{F_\Th}\;\;.
\end{eqnarray*}
Due to part 1), every element of $\Spm\FK{G}$ is of the form $g\tr e(R(\Xi))h$ with $g,h\in G_f$, $\Xi$ special. The equations
\begin{eqnarray*}
  0\;\,=\;\, f_{vw}(g\tr e(R(\Xi))h)\;\,=\;\,\kB{gv}{e(R(\Xi))hw} &\mb{ for all }& v,w\in L(\La)\,,\;\La\in P^+\setminus 
  \overline{F_\Th}
\end{eqnarray*}
are equivalent to
\begin{eqnarray*}
    e(R(\Xi))L(\La)\;\,=\;\,\{0\} &\mb{ for all }& \;\La\;\in\; P^+\setminus \overline{F_\Th}\;=\;P^+\setminus R(\Th)\;\;.
\end{eqnarray*}
Recall that for an element $\La\in P^+$ we have $P(\La)\cap R(\Xi)=\emptyset$ if and only if $\La\notin R(\Xi)$. Therefore these equations are 
equivalent to $\La\in P^+\setminus \overline{F_\Xi}$ for all $\La\in P^+\setminus \overline{F_\Th}$. This is equivalent to 
$P^+\setminus \overline{F_\Th}\subseteq  P^+\setminus \overline{F_\Xi}$,  which in turn is equivalent to $\Th\subseteq \Xi$.\vspace*{1ex}\\
c) The closure of the $G_f\times G_f$-orbit $G_f\tr e(R(\Th))G_f$ is a union of $G_f\times G_f$-orbits. Due to b) it is 
contained in $\bigcup_{\Xi\;sp.\;,\;\Xi\supseteq \Th} G_f\tr e(R(\Xi))G_f$. To show equality, it is sufficient to show, that 
the closure contains the elements $1\tr e(R(\Xi))$, $\Xi\supseteq\Th$, $\Xi$ special.\\
Because left multiplications with elements of $\GD$ are Zariski continuous on $\GD$, we find
\begin{eqnarray*}
   e(R(\Th))\GD \;\,=\;\, e(R(\Th))\Aa{G}\;\,\subseteq\;\, \Aa{e(R(\Th))G}\;\;.
\end{eqnarray*}
Therefore we get
\begin{eqnarray*}
   1\tr e(R(\Th))\GD \;\,\subseteq\;\, \Asp{1\tr e(R(\Th))G_f}\;\,\subseteq\;\,\Asp{G_f\tr e(R(\Th))G_f}\;\;.
\end{eqnarray*}
Now $e(R(\Th))\GD$ contains for every $\Xi\supseteq \Th$, $\Xi$ special, the element $e(R(\Th))e(R(\Xi))=
e(R(\Th)\cap R(\Xi))=e(R(\Xi))$.\vspace*{1ex}\\ 
d) If $\Th$ is the biggest special set with respect to the inclusion, then the orbit $G_f\tr e(R(\Th))G_f$ is closed due to c).\\
If $\Th$ is not the biggest special set, then due to c) we find
\begin{eqnarray*}
   G_f\tr e(R(\Th))G_f &=& \Asp{ G_f\tr e(R(\Th))G_f}\,\setminus\,\bigcup_{\Xi\,sp.\,,\,\Xi\supsetneqq \Th}\Asp{G_f\tr e(R(\Xi))G_f}\;\;.
\end{eqnarray*}
There are only finitely many special sets. Therefore $G_f\tr e(R(\Th))G_f$ is locally closed.\vspace*{1ex}\\
e) The algebra of strongly regular functions $\FK{G}$ is an integral domain. Therefore every subset of $\Ae{G}=\widehat{G_f}$, which contains $G$, 
is irreducible. In particular $G_f$ is irreducible. From Theorem \ref{SettingirreducibleOrbits} follows, that the $G_f\times G_f$-orbits of 
$\Spm\FK{G}$ are irreducible.\\
\End
Our next aim is to define and describe big cells of every $G_f\times G_f$-orbit of $\Spm\FK{G}$. As a preparation we first determine certain 
principal open sets of $\Spm\FK{G}$.\\ 
Recall that $\gt_\La$ denotes the function defined by $\gt_\La:=\frac{1}{\kB{v_\La}{v_\La}} f_{v_\La v_\La}$, 
where $v_\La\in L(\La)_\La\setminus\{0\}$, $\La\in P^+$. Recall the multiplicative property $\gt_\La\gt_{\La'}=\gt_{\La+\La'}$, $\La,\La'\in P^+$.
\begin{Theorem}\label{S1b} Let $\Th$ be special. The principal open subset of $\Spm\FK{G}$ associated with $\gt_\La$, $\La\in F_\Th\cap P$, is given by 
\begin{eqnarray}\label{principalopen}
     \dot{ \bigcup_{\hat{\sigma}\in \widehat{\cal W}_\Th} } U_f\tr (\hat{\sigma} T) U_f \;\;.
\end{eqnarray}
This set, as well as its coordinate ring as a principal open set, is independent of the chosen element $\La\in F_\Th\cap P$.
\end{Theorem}
We denote this principal open set by $D(\Th)$, and its coordinate ring as a principal open set by $\FK{D(\Th)}$.\vspace*{1ex}\\
\Proof
a) We first show that the principal open set of $\gt_\La$ is given by (\ref{principalopen}):
Due to the Birkhoff decomposition, every element of $\Spm\FK{G}$ can be written in the form 
$u\tr n_\sigma e(R)\ti{u}$, with $u,\ti{u}\in U_f$, $n_\sigma\in N$, and $R$ a face of $X$.
Let $v_\La\in L(\La)_\La\setminus\{0\}$. We find
\begin{eqnarray*}
   0\;\,\neq\;\, \gt_\La (u\tr n_\sigma e(R)\ti{u})\;\,=\;\,\frac{ \kB{uv_\La}{n_\sigma e(R) \ti{u} v_\La}   }{\kB{v_\La}{ v_\La}}
  \;\,=\;\, \frac{ \kB{v_\La}{n_\sigma e(R) v_\La}  }{\kB{v_\La}{ v_\La}}
\end{eqnarray*}
if and only if $\La\in R$ and $\sigma \La=\La$. Because $\La$ is an interior point of the face $R(\Th)$, the first condition is equivalent to 
$R(\Th)\subseteq R$. The second condition is equivalent to $\sigma\in\We_\Th$.\\ 
Because of
\begin{eqnarray*}
  \WeD_\Th &=& \dot{\bigcup_{S\;a\; face\; of\; X\atop S\supseteq R(\Th)}} \We_\Th \;\ve{S}\;\;,
\end{eqnarray*}
we get $\sigma \ve{R}\in\WeD_\Th$. We also find that for every $\hat{\sigma}\in\WeD_\Th$ we have 
$U_f\tr \hat{\sigma}T U_f\subseteq D(\gt_\La)$.\vspace*{1ex}\\
b) Clearly (\ref{principalopen}) does not depend on $\La\in F_\Th\cap P$.  
To show that the coordinate rings of the principal open sets are independent of $\La\in F_\Th\cap P$, we only have to show, that for any 
$N, N'\in F_\Th\cap P$ there exists a function $f$ of the coordinate ring of $\Spm\FK{G}$, and an integer $n\in\N$, such that
\begin{eqnarray}\label{fquotient}
     \frac{f}{(\gt_N)^n}\res{D(\Th)}&=& \frac{1}{\gt_{N'}}\res{D(\Th)}\;\;.
\end{eqnarray}
Because $F_\Th$ is open in the linear span of $F_\Th$, there exists an integer $n\in\N$, such that $N-\frac{1}{n}N'\in F_\Th$. Since 
$F_\Th$ is a cone, we find $nN-N'\in F_\Th\cap P$. The function $f=\gt_{nN-N'}$ satisfies equation (\ref{fquotient}).\\
\End
Let $\Th$ be special. We call the set $BC(\Th):=U_f \tr e(R(\Th))T U_f$, as well as every translate $U_f g\tr e(R(\Th))T U_f h$, where $g,h\in G_f$, a 
big cell of the orbit $G_f \tr e(R(\Th))G_f$. This name is justified by the following three theorems.
\begin{Theorem} \label{S2} Let $\Th$ be special. The big cell $BC(\Th)$ is principal open in the closure of $G_f\tr e(R(\Th))G_f$, i.e., 
\begin{eqnarray}\label{bcint} 
     BC(\Th) &=& \Asp{G_f\tr e(R(\Th))G_f}\cap D(\Th)\;\;.
\end{eqnarray}
The big cell $BC(\Th)$ is dense in the closure of $G_f\tr e(R(\Th))G_f$.
\end{Theorem}
\Proof
a) We first show formula (\ref{bcint}). Due to the Birkhoff decomposition of $\Spm\FK{G}$, and the formulas (\ref{closure}), (\ref{orbit}), and 
(\ref{principalopen}), we only have to show  
\begin{eqnarray*}
   \{\,\ve{R(\Th)}\,\}  &=&  \bigcup_{\Xi\;sp.\,,\,\Xi\supseteq \Th}\We\,\ve{R(\Xi)} \,\We \;\;\cap\;\; \WeD_\Th  \;\;. 
\end{eqnarray*}
Because of $\WeD_\Th = \bigcup_{\Xi\,sp.,\,\Xi\subseteq\Th} \We_\Th\,\ve{R(\Xi)}\, \We_\Th$, the intersection on the right equals 
$\We_\Th\,\ve{R(\Th)}\,\We_\Th$. Since $\We_\Th$ is the pointwise stabilizer of $R(\Th)$, its elements fix $\ve{R(\Th)}$. Therefore this intersection 
contains only the element $\ve{R(\Th)}$.\vspace*{1ex}\\ 
b) Due to Theorem \ref{S1} the orbit $G_f\tr e(R(\Th))G_f$ is irreducible. Therefore also its closure is irreducible. Because the big cell is open 
in the closure, and nonempty, it is dense.\\
\End
Let $\Th$ be special. We equip the big cell $BC(\Th)$ with its coordinate ring $\FK{BC(\Th)}$ as a principal open set in the closure of 
$G_f\tr e(R(\Th))G_f$.\\
Set $T^\Th:=T_{I\setminus \Th}T_{rest}$. Set $P^\Th:=\Z\mb{-}span\Mklz{\La_i}{i=1,\ldots,2n-l,\;i\notin \Th }$, 
and identify the group algebra $\FK{P^\Th}$ with the classical coordinate ring of the torus $T^\Th$.\\ 
Note that due to Theorem \ref{SpU2} the coordinate ring $\FK{U_f^\Th}$ is isomorphic to $\FK{U^\Th}$ by the restriction map.
\begin{Theorem}\label{S3} Let $\Th$ be special. We get an isomorphism
\begin{eqnarray*}
  m:\;U_f^\Th \times T^\Th \times U_f^\Th &\to & BC(\Th)
\end{eqnarray*}
by $m(u,t,\ti{u}):=u\tr e(R(\Th)) t\ti{u}$, where $u,\,\ti{u}\in U_f^\Th$ and $t\in T^\Th$.
\end{Theorem}
\Proof
a) First we show, that $m$ is surjective:
Due to \cite{M}, Proposition 2.13, compare the section preliminaries, we have $e(R(\Th))T_\Th=e(R(\Th))$. Because of $T=T_\Th T^\Th$ we get 
$e(R(\Th))T=e(R(\Th))T^\Th$.\\
Due to the same proposition we also have $e(R(\Th))U_\Th=e(R(\Th))$. Left and right multiplications with elements of 
$\Aa{G}=\GD$ are Zariski continuous on $\Ae{G}=\widehat{G_f}$. Therefore we get by using Theorem \ref{SpU2}:
\begin{eqnarray*}
 e(R(\Th)) (U_f)_\Th\;\,=\;\, e(R(\Th))\Ae{U_\Th} \;\,\subseteq \;\, \Ae{e(R(\Th))U_\Th}\;\,=\;\, 
\Ae{\{e(R(\Th))\} } \;\,=\;\, \{e(R(\Th))\}\;\;. 
\end{eqnarray*}
Because of $U_f=(U_f)_\Th\ltimes (U_f)^\Th$ we find $e(R(\Th))U_f=e(R(\Th))(U_f)^\Th$.\\
Because of these formulas, and because $\TD$ is abelian, we get:
\begin{eqnarray*}
  U_f\tr e(R(\Th))TU_f &=& U_f\tr e(R(\Th)) T^\Th (U_f)^\Th \;\,= \;\;e(R(\Th))U_f\tr T^\Th (U_f)^\Th \\ 
                       &=& (U_f)^\Th\tr e(R(\Th)) T^\Th (U_f)^\Th\;\;.
\end{eqnarray*}
b) Next we show that the comorphism $m^*: \FK{BC(\Th)} \to \FK{U_f^\Th}\otimes \FK{P^\Th}\otimes \FK{U_f^\Th}$ is well defined and surjective: 
For $\La\in P^+$ choose $\kBl$-dual bases of $L(\La)$, by choosing $\kBl$-dual bases 
\begin{eqnarray*}
  (a_{\la i})_{i=1,\ldots,m_\la} &\quad,\quad & (b_{\la i})_{i=1,\ldots,m_\la}
\end{eqnarray*}
of $L(\La)_\la$ for every $\la\in P(\La)$. Fix elements $N\in F_\Th\cap P$ and $v_N\in L(N)_N\setminus\{0\}$. For $v,w\in L(\La)$, 
and $u,\ti{u}\in U_f^\Th$, $t\in T^\Th$ we get
\begin{eqnarray*}
\frac{f_{vw}}{(\gt_N)^k}\,(u\tr e(R(\Th))t\ti{u})\;\,=\;\, 
  \sum_{\la\in P(\La)\cap R(\Th),\;i} \kB{uv}{a_{\la i}}e_{\la-k N}(t) \kB{b_{\la i}}{\ti{u} w}\;\;.
\end{eqnarray*}
This sum has only finitely many nonzero summands, because $supp(uv)$ and $supp(\ti{u}w)$ are finite. Denote by $p^\Th:P\to P^\Th$ the projection 
corresponding to the decomposition $P=P_\Th\oplus P^\Th$. Due to the last formula we have
\begin{eqnarray*}
\frac{f_{vw}}{(\gt_N)^k}    \res{BC(\Th)}\,\circ\, m\;\,=\;\, 
  \sum_{\la\in P(\La)\cap R(\Th),\;i}  f_{a_{\la i} v}\res{U^\Th}\otimes\, e_{p^\Th(\la)-kN}\otimes f_{b_{\la i} w}\res{U^\Th}\;\;.
\end{eqnarray*}
There are only finitely many nonzero summands of this sum. The function $f_{a_{\la i}v}\res{U^\Th}$ is nonzero at most if $v\neq 0$ and 
if $\la$ is bigger than a weight of $supp(v)$, which is only possible for finitely many weights $\la$ in $P(\La)$. Similar things hold for 
$f_{b_{\la i}w}\res{U^\Th}$.\\ 
In particular, $m^*$ is well defined. From this formula we find for $v\in L(N)$ and $\La\in\overline{F_\Th}\cap P$:
\begin{eqnarray*}
   m^*(\frac{f_{v v_N}}{\gt_N}) &=& f_{v_N v}\res{U^\Th}\otimes\, 1\otimes 1 \;\;,\\
   m^*(\frac{\gt_\La}{(\gt_N)^k}) &=&  1\otimes e_{\La-k N}\otimes 1 \;\;,\\
   m^*(\frac{f_{v_N v}}{\gt_N}) &=& 1\otimes 1\otimes f_{v_N v}\res{U^\Th}\;\;.
\end{eqnarray*}
It is easy to see, that $(\overline{F_\Th}\cap P)-\Nn N= P^\Th$. Taking into account Theorem \ref{SpU1}, 2), we have found elements of 
$\FK{BC(\Th)}$, which are mapped onto a system of generators of $\FK{U_f^\Th}\otimes \FK{P^\Th}\otimes \FK{U_f^\Th}$. Therefore
$m^*$ is surjective.\vspace*{1ex}\\
c) $m^*$ is injective, because $m$ is surjective. To show the injectivity of $m$, let $u_1,u_2,\ti{u}_1,\ti{u}_2\in U_f^\Th$, 
$t_1,t_2\in T^\Th$ such that $m(u_1,t_1,\ti{u}_1)=m(u_2,t_2,\ti{u}_2)$. Then for all $f\in \FK{U^\Th}$ we have
\begin{eqnarray*}
  f(u_1)\;\,=\;\, (m^*)^{-1}(f\otimes 1\otimes 1)(\underbrace{m(u_1,t_1,\ti{u}_1)}_{=m(u_2,t_2,\ti{u}_2)}) \;\,=\;\, f(u_2)\;\;.
\end{eqnarray*}
Therfore we find $u_1=u_2$. In a similar way, we get $t_1=t_2$, and $\ti{u}_1=\ti{u}_2$.\\
\End
In the next theorem we give countable coverings of every $G_f\times G_f$-orbit of $\Spm\FK{G}$ by big cells.\\
Recall that $\We^J$ denotes the set of minimal coset representatives of $\We/\We_J$, and $\mb{}^J\We$ denotes the set of minimal coset representatives of 
$\We_J\backslash \We$, $J\subseteq I$.
\begin{Theorem} Let $\Th$ be special. We have
\begin{eqnarray*}
  G_f \tr e(R(\Th))G_f  &=&  \bigcup_{\sigma\in \mb{}^\Th{\cal W}\;,\;\tau\in \mb{}^{\Th\cup\Th^\bot}{\cal W}} U_f
  \sigma \tr e(R(\Th))T  U_f \tau \\
                    &=&   \bigcup_{\sigma\in \mb{}^{\Th\cup\Th^\bot}{\cal W}\;,\;\tau\in \mb{}^\Th{\cal W} } U_f 
  \sigma \tr e(R(\Th))T U_f \tau     \;\;.
\end{eqnarray*}
\end{Theorem}
\Proof
It is sufficient to prove the second covering of the theorem, the first follows by applying the Chevalley involution of $\Spm\FK{G}$.
Obviously the sets $U_f\sigma \tr e(R(\Th))T U_f \tau $ are contained in the orbit $G_f\tr e(R(\Th))G_f$. Therefore it is sufficient to show that 
the orbit $G_f\tr e(R(\Th)) G_f$ is contained in the second union.\\
By writing the elements $\We\ve{R(\Th)}\We$ in normal form, and inserting in (\ref{orbit}), we get
\begin{eqnarray*}
   G_f\tr e(R(\Th))G_f &=& \bigcup_{    \sigma\in{\cal W}^{\Th\cup\Th^\bot},\;\tau\in \mb{}^\Th{\cal W}    } 
   U_f\tr \sigma e(R(\Th))T\tau U_f\;\;.
\end{eqnarray*}
To transform the expression $U_f\tr \sigma e(R(\Th))T\tau U_f$, we use the following decomposition of $U_f$, associated to an 
element $w\in\We$, proved in \cite{Sl1}, Section 5.5:\\ 
Set $\W^+_w:=\Mklz{\al\in \pW}{ w\al\in\nW}$ and $(\pW)^w:=\Mklz{\al\in \pW}{ w\al\in\pW}$. Then $\W^+_w$ consists of
finitely many positive real roots, and we have
\begin{eqnarray*}
 U_f \;\,=\;\,U_w\,U_f^w &\mb{ where }& U_w:=\exp(\bigoplus_{\al\in \Delta_w^+}\g_\al) \;,\;
   U_f^w:=\exp(\prod_{\al\in (\Delta^+)^w}\g_\al)\;\;.
\end{eqnarray*} 
Using this decomposition, and Proposition 2.14 of \cite{M}, compare the section preliminaries, we find
\begin{eqnarray*}
  U_f\tr \sigma e(R(\Th)) T \tau U_f \;\,=\;\, U_f\tr\sigma e(R(\Th)) T\underbrace{\tau U_\tau \tau^{-1}}_{\subseteq U^-} 
\underbrace{\tau U_f^\tau\tau^{-1}}_{\subseteq U_f}\tau\\
\subseteq\;\,U_f\tr \sigma U^-_{\Th^\bot} e(R(\Th)) T U_f\tau \;\,=\;\, 
  (\sigma U_{\Th^\bot}^-\sigma^{-1})^* U_f\tr \sigma e(R(\Th)) T U_f\tau\;\;. 
\end{eqnarray*}
Because of $\sigma\in \We^{\Th\cup\Th^\bot}=\Mklz{\sigma\in\We}{\sigma\al_i\in\prW \mb{ for all } i\in \Th\cup\Th^\bot}$, we have 
$\sigma U^-_{\Th^\bot}\sigma^{-1}\subseteq U^-$, and the last expression equals $U_f\tr \sigma e(R(\Th))T U_f\tau$. By doing similar 
transformations, we get
\begin{eqnarray*}
 \lefteqn{ U_f\tr \sigma e(R(\Th))T U_f\tau \;\,=\;\,\sigma^{-1}U_f\tr e(R(\Th))T U_f\tau } \\
  &=& \underbrace{(\sigma^{-1}U_{\sigma^{-1}}\sigma)}_{\subseteq U^-} 
             \underbrace{(\sigma^{-1} U_f^{\sigma^{-1}}\sigma )}_{\subseteq U_f}\sigma^{-1} \tr e(R(\Th))T U_f\tau 
  \;\,\subseteq\;\, U_f\sigma^{-1}\tr U e(R(\Th)) T U_f\tau \\
   &=& U_f\sigma^{-1}\tr e(R(\Th))T U_{\Th^\bot}U_f\tau 
  \;\,=\;\, U_f\sigma^{-1}\tr e(R(\Th))T U_f\tau \;\;.
\end{eqnarray*}
Because of $(\We^{\Th\cup\Th^\bot})^{-1}=\mb{}^{\Th\cup\Th^\bot}\We$ we have shown, that the orbit is contained in the second union.\\
\End
Our last aim is to show, that there exist stratified transversal slices to the $G_f\times G_f$-orbits of $\Spm\FK{G}$ at any of their points.\\
Because $G_f\times G_f$ acts by isomorphisms on $\Spm\FK{G}$, it is sufficient to find stratified transversal slices at the points 
$1\tr e(R(\Th))\in G_f\tr e(R(\Th)) G_f$, $\Th$ special.\vspace*{1ex}\\ 
As transversal slice at $1\tr e(R(\Th))$, we will use the closure $\Asp{G_\Th}:= \Asp{1\tr G_\Th}$, equipped with its coordinate 
ring as a closed subset of $\Spm\FK{G}$. We have the following description:
\begin{Theorem} Let $\Th$ be special. The restriction map $\FK{G}\to\FK{G_\Th}$ induces a closed embedding 
$\Spm\FK{G_\Th}\to\Spm\FK{G}$ with image
\begin{eqnarray*}
   \Asp{G_\Th} &=& (\widehat{G_f})_\Th\tr (\widehat{G_f})_\Th\;\,=\;\,
                       \dot{\bigcup_{\Xi\subseteq \Th,\,\Xi\,special}} (G_f)_\Th\tr e(R(\Xi)) (G_f)_\Th   \\
               &=&   \dot{ \bigcup_{\hat{\sigma}\in \widehat{\cal W}_\Th} } (U_f)_\Th\tr (\hat{\sigma} T_\Th) (U_f)_\Th \;\;.
\end{eqnarray*} 
\end{Theorem}
\Proof
It is not difficult to check, that the map $\Spm\FK{G_\Th}\to\Spm\FK{G}$ is a closed embedding  with image $\Asp{G_\Th}$. 
Due to Theorem \ref{SpG1}, $(\widehat{G_f})_\Th\times (\widehat{G_f})_\Th$ maps surjectively to $\Spm\FK{G_\Th}$. To show 
the first equation of the theorem, we have to show, that the concatenation of the maps 
\begin{eqnarray*}
    (\widehat{G_f})_\Th\times (\widehat{G_f})_\Th \;\,\stackrel{\al}{\twoheadrightarrow}\;\, \Spm\FK{G_\Th} \;\,\to\;\, \Spm \FK{G}
\end{eqnarray*}
coincides with the restricted map $\tr :(\widehat{G_f})_\Th\times (\widehat{G_f})_\Th\to \Spm\FK{G}$.\\
Let $x,y\in (\widehat{G_f})_\Th$. Let $v,w\in L(\La)$, $\La\in P^+$. Choose a decomposition $L(\La)=\bigoplus_{j\in J} V_j$ of $L(\La)$ in an 
orthogonal direct sum of irreducible highest weight $(\g_\Th+\h)$-modules. Write $v$, $w$ as sums $v=\sum_{j\in J} v_j$, $w=\sum_{j\in J} w_j$ 
with $v_j,w_j\in V_j$. Then:
\begin{eqnarray*}
 \al(x,y)(f_{vw}\res{G_\Th}) &=& \al(x,y)(\,\sum_{j\in J}f_{v_j w_j}\res{G_\Th}\,)\;\,=\;\,\sum_{j\in J} \kB{xv_j}{yw_j}\\
                             &=&\kB{xv}{yw}\;\,=\;\, (x\tr y) (f_{vw})
\end{eqnarray*}
The other equations of the theorem follow from
\begin{eqnarray*}
&&(\widehat{G_f})_\Th\tr(\widehat{G_f})_\Th \;=\; U^-_\Th \ND_\Th (U_f)_\Th \tr (\widehat{G_f})_\Th \;=\; 
      (U_f)_\Th \tr (U^-_\Th \ND_\Th )^* (\widehat{G_f})_\Th \\ 
&&  =\;  (U_f)_\Th \tr (\widehat{G_f})_\Th \;=\; (U_f)_\Th \tr  U^-_\Th \ND_\Th (U_f)_\Th \;=\; (U^-_\Th)^*(U_f)_\Th \tr \ND_\Th (U_f)_\Th \\ 
&& =\;(U_f)_\Th \tr \ND_\Th (U_f)_\Th\; \subseteq \;\bigcup_{\Xi\subseteq \Th,\,\Xi\,special} (G_f)_\Th\tr e(R(\Xi)) (G_f)_\Th \;\subseteq\;  
     (\widehat{G_f})_\Th\tr(\widehat{G_f})_\Th\;\;. 
\end{eqnarray*}
The unions in the equations of the theorem are disjoint, because the unions $\bigcup_{\Xi\, special} G_f\tr e(R(\Xi)) G_f$, 
$\bigcup_{\hat{w}\in\hat{\cal W}} U_f\tr \hat{w}T U_f$ are disjoint.\\
\End
As open neighborhood of $1\tr e(R(\Th))$, we take the principal open set $D(\Th)$.
\begin{Theorem} Let $\Th$ be special.\\
1) We have $G_f\tr e(R(\Th))G_f\cap D(\Th)=BC(\Th)$. The big cell $BC(\Th)$ is closed in the principal open set $D(\Th)$. Its coordinate ring 
coincides with the coordinate ring as a closed subset of the principal open set $D(\Th)$.\vspace*{1ex}\\ 
2) We have $1\tr e(R(\Th))\in\Asp{G_\Th}\subseteq D(\Th)$. The coordinate ring of $\Asp{G_\Th}$ coincides with the coordinate ring as a closed 
subset of the principal open set $D(\Th)$.\vspace*{1ex}\\
3 a) We get an isomorphism 
\begin{eqnarray*}
 \Psi:\; \Asp{G_\Th}\,\times \,BC(\Th) &\to & D(\Th)
\end{eqnarray*}
by $\Psi(x\tr y, u\tr  e(R(\Th))t\ti{u}):= xu\tr yt\ti{u}$, where $x,y\in(\widehat{G_f})_\Th$, $u,\ti{u}\in (U_f)^\Th$, and $t\in T^\Th$.\vspace*{1ex}\\
b) Inserting $1\tr e(R(\Th))$ in the first entry (resp. second entry) of $\Psi$ induces the identity map on $BC(\Th)$ (resp. $\Asp{G_\Th}$).\\
c) The partition of $\Spm\FK{G}$ into $G_f\times G_f$-orbits induces partitions of $\Asp{G_\Th}$ and $D(\Th)$. $\Psi$ preserves the orbits, i.e., 
\begin{eqnarray*}
   \Psi(\,\Asp{G_\Th}\cap G_f\tr e(R(\Xi))G_f\,,\, BC(\Th)\,) \;\,=\;\,D(\Th)\,\cap\,  G_f\tr e(R(\Xi))G_f \quad,\quad \Xi \mb{ special}\;\;.
\end{eqnarray*}
\end{Theorem}
\Proof 
1) The first part of the theorem follows immediately from the definition of the big cell and its coordinate ring.\vspace*{1ex}\\
2) From the last theorem and the description of $D(\Th)$ given in Theorem \ref{S1b} follows $1\tr e(R(\Th))\in\Asp{G_\Th}\subseteq D(\Th)$.\\
The statement of part 2) about the coordinate rings can be seen as follows:
$D(\Th)$ is the principal open set of $\gt_\La$ for an element $\La\in F_\Th \cap P$. Now $G_\Th$ stabilizes every point of the highest weight space 
$L(\La)_\La$, because for $j\in \Th$ the $\al_j$-string of $P(\La)$ through $\La$ consists only of $\La$. Therefore $\gt_\La$ takes the 
value 1 on $1\tr G_\Th$, and also on the closure $\Asp{1\tr G_\Th}=\Asp{G_\Th}$.\vspace*{1ex}\\
3) Because of Theorem \ref{S3}, the statement of 3 a) is equivalent to the following statement, which we will prove: We get a bijective map 
\begin{eqnarray*}
 \ti{\Psi}:\, (U_f)^\Th\times \Asp{G_\Th}\times T^\Th \times (U_f)^\Th  &\to & D(\Th)
\end{eqnarray*}
by $\Psi(u, x\tr y, t, \ti{u}):= xu\tr yt\ti{u}=(x\tr y)(u,t\ti{u})$, where $x,y\in(\widehat{G_f})_\Th$, $u,\ti{u}\in (U_f)^\Th$, $t\in T^\Th$. Its 
comorphism exists and is an isomorphism of algebras.\vspace*{1ex}\\
From the the descriptions of $D(\Th)$ and $\Asp{G_\Th}$ given in Theorem  \ref{S1b} and in the last theorem follows, that the image of the map 
$\ti{\Psi}$ is $D(\Th)$.\vspace*{1ex}\\
Next we show that the comorphism $\ti{\Psi}^*$ exists and is surjective: Choose an element $\La\in F_\Th\cap P_I$. It is easy to check, that the 
multiplication maps $D_{G'}(\gt_\La)\times T_{rest}\to D_G(\gt_\La)$, and 
$T_{I\setminus\Th}\times T_{rest}\to T^\Th$ are bijective. Furthermore their comorphisms are isomorphisms 
$\FK{D_G(\gt_\La)}\to \FK{D_{G'}(\gt_\La)}\otimes \FK{P_{rest}}$, and $\FK{P^\Th}\to \FK{P_{I\setminus\Th}}\otimes \FK{P_{rest}}$. Also the comorphism of the 
bijective map $*:U^\Th\to (U^\Th)^-$ is an isomorphism $\FK{(U^\Th)^-}\to \FK{U^\Th}$. Taking into account Theorem \ref{Ueb2} for $J=I$ and $L=\Th$, 
we find that the map 
\begin{eqnarray*}
  \ti{m}:\, U^\Th\times G_\Th \times T^\Th \times U^\Th &\to & D_G(\gt_\La)
\end{eqnarray*}
given by $\ti{m}(u,g,t,\ti{u}):=u^* g t \ti{u}$ is bijective, its comorphism exists, and is an isomorphism of algebras.\\
Now identify $G$ with $1\tr G$. Then $D_G(\gt_\La)$ is contained in $D(\Th)$. Due to Theorem \ref{SpU2} the set $U^\Th$ is dense in 
$(U_f)^\Th$. The coordinate ring $\FK{(U_f)^\Th}$ is isomorphic to $\FK{U^\Th}$ by the restriction map. Similar things hold for 
$G_\Th$, $\Asp{G_\Th}$ and their coordinate rings.\\
For a coordinate ring $\FK{B}$ and a nonempty subset $A\subseteq B$ denote by $res_A^B$ the restriction map $\FK{B}\to \FK{A}$. It is easy to check, 
that the surjective map
\begin{eqnarray*}
 \left(  res_{\,U^\Th}^{(U_f)^\Th}\,\otimes\, res_{G_\Th}^{\Asp{G_\Th}} \,\otimes\,res_{T^\Th}^{T^\Th}\,\otimes\,res_{\,U^\Th}^{(U_f)^\Th}\right)^{-1}
    \;\circ\;\ti{m}^*\;\circ\;res_{D_G(\gt_\La)}^{\,D(\Th)} 
\end{eqnarray*}
is the comorphism of $\ti{\Psi}$. (Make use of $\ti{\Psi}(u,g,t,\ti{u})=\ti{m}(u,g,t,\ti{u})$ for $u,\ti{u}\in U^\Th$, $g\in G_\Th$, and $t\in T^\Th$.)
\vspace*{1ex}\\
Because the maps $\ti{\Psi}$ and $\ti{\Psi}^*$ are surjective, they are also injective. This is shown in the same way as the injectivity of the maps 
$m$ and $m^*$ in the proof of Theorem \ref{S3}.\vspace*{1ex}\\ 
3 b) follows immediately from the definition of $\Psi$. 3 c) can be checked easily by using the definition and the bijectivity of $\Psi$, and the 
last theorem.\\ 
\Ende
%
%
%
%
%
\mb{}\\
{\bf Acknowledgment} I would like to thank the Deutsche Forschungsgemeinschaft for providing my main financial support, when most of 
this paper has been written. I also would like to thank the TMR-Program ERB FMRX-CT97-0100 ``Algebraic Lie Theory'' for providing some financial 
support for traveling.\\
Furthermore I would like to thank the referees for their useful comments. 
%
%
%
%
%
%
%

%
%
%
\end{document}